%% file: SM-ROM.tex
\documentclass[11pt,feqno]{article}
\usepackage{amssymb,bbm,amsmath,amsthm,amscd}
\usepackage{vmargin}
\usepackage[dvipsnames]{xcolor}
\usepackage{array}
\usepackage{soul}
 
 \usepackage{graphicx}
\DeclareGraphicsExtensions{.pdf,.eps}

\def\cb{\color{black}}
\def\cn{\color{black}}
\catcode`\@=11
\newbox\bz@
\newdimen\bdimz@
\def\linethrough#1{\setbox\bz@=\hbox{#1}
\bdimz@=\ht\bz@ \divide\bdimz@ by 5 \advance\bdimz@ by -\dp\bz@ \ht\bz@=\bdimz@
\leavevmode\hbox{$\overline{\overline{\box\bz@}}$\relax}}
\catcode`\@=12

\def\nbZ{{\mathchoice {\hbox{$\sf\textstyle Z\kern-0.4em Z$}}
{\hbox{$\sf\textstyle Z\kern-0.4em Z$}}
{\hbox{$\sf\scriptstyle Z\kern-0.3em Z$}}
{\hbox{$\sf\scriptscriptstyle Z\kern-0.2em Z$}}}}

\newcounter{compteur}[subsection]

\usepackage[colorlinks=true]{hyperref}

\newtheorem{theorem}{Theorem}[section]
\newtheorem{lemma}[theorem]{Lemma}

\newtheorem{remark}[theorem]{Remark}

\newcommand{\bu}{{\boldsymbol u}}

\usepackage{macros}

\makeatletter 

\@addtoreset{equation}{section}

\usepackage{tikz}
\usetikzlibrary{positioning}
\usepackage{caption}
\usepackage{subcaption}
\usepackage{float}
\usepackage{bm}
\usepackage[export]{adjustbox}
\input{./NEW_COMMANDS.tex}
\newcommand{\ModifMourad}[1]{\textcolor{black}{#1}}
\newcommand{\review}[1]{\textcolor{black}{#1}}
\begin{document}


\author{
Tom\'as Chac\'on Rebollo \thanks{Departamento EDAN \& IMUS, Universidad de Sevilla, Spain. {\tt chacon@us.es}}, 
Samuele Rubino \thanks{Departamento EDAN \& IMUS, Universidad de Sevilla, Spain. {\tt samuele@us.es}},
Mourad Oulghelou \thanks{LAMPA, Arts et M\'etiers ParisTech, France. {\tt mourad.oulghelou@ensam.eu}},
Cyrille Allery \thanks{LaSIE, Universit\'e de la Rochelle, France. {\tt cyrille.allery@univ-lr.fr}}}

\title{Error analysis of a residual-based stabilization-motivated POD-ROM for incompressible flows}
\maketitle

\begin{abstract}
\ModifMourad{This article presents error bounds for a velocity-pressure segregated POD reduced order model discretization of the Navier-Stokes equations. The stability is proven in $L^\infty(L^2)$ and energy norms for velocity, with bounds that do not depend on the viscosity, while for pressure it is proven in a semi-norm of the same asymptotic order as the $L^2$ norm with respect to the mesh size. The proposed estimates are calculated for the two flow problems, the flow past a cylinder and the lid-driven cavity flow. Their quality is then assessed in terms of the predicted logarithmic slope with respect to the velocity POD contribution ratio. We show that the proposed error estimates allow a good approximation of the real errors slopes and thus a good prediction of their rate of convergence. }
\end{abstract}


\medskip


\section{Introduction}\label{sec:Intro}
\ModifMourad{Since the beginning of the XXI Century, there is an increasing number of applications of Reduced Order Models (ROMs) addressed to the approximation of incompressible Navier-Stokes flow solutions. By using carefully chosen basis functions forming the projection subspace, ROMs \review{collaborating with high fidelity solutions}, offer a computationally efficient \review{tool to reduce the cost of} heavy numerical \review{simulations}.
A very famous method that is extensively used to construct projection subspaces is the POD (Proper Orthogonal Decomposition) method \cite{Sirovich, Berkooz}. Starting from a set of "snapshots" computed from the high-fidelity model (or full-order model, FOM) in the off-line phase, the POD allows trough solving an eigenvalue problem to obtain a truncated basis formed by modes corresponding to the first larger eigenvalues. Once the POD basis is built, the temporal dynamics of the targeted system (in our case, the Navier-Stokes equations) is determined by solving a low dimension system obtained from the Galerkin projection of the model equations onto this reduced subspace \cite{Rowley}. 
The resulting model provides very fast solutions of the targeted equations, with computational cost frequently several order of magnitude lower than that of the FOM. The POD Galerkin approach has been applied in many fields such as aeroelasticity \cite{Thomas2003, Lucia2004}, structural dynamics \cite{Park2007,ALDMOUR2002,Buljak2011}, modal analysis \cite{Han2003}, flow prediction and control \cite{Allery-2005, Beghein-Allery-2014, Fick2018, Hijazi2020, Tallet2016}, particles dispersion \cite{Allery2008, Allery2005, Beghein2014}, etc. Other approaches to construct projection subspaces relying on greedy algorithms can be found in \cite{Ballarin2020,AlfioRozza2009}. For a thorough review, the interested reader can refer to \cite{Hesthaven2016}.
%
%
%
%
%
%
%
%
%
%
%
\\
\noindent In practice, the velocity snapshots are weakly divergence free up to the computer accuracy. This allows to drop the pressure from projected equations and form a reduced order model that only contains the velocity amplitudes as unknown variables. However for many applications, the pressure is needed to be recovered. Namely, to calculate lift and drag forces exerted on walls or immersed boundaries. One option to recover the pressure is to build a ROM that also accounts for its time trajectory in a reduced coupled velocity-pressure setting. Specific techniques to ensure the discretisation of the pressure are thus needed, for instance using velocity sumpremizers of the velocity-pressure duality \cite{Quarter2015, Veroy2007} or stabilization techniques \cite{Rubino20}. 
An alternative way is to recover the pressure once the velocity is computed from the velocity ROM. The governing reduced order model in this case is a pressure Poisson equation derived by taking the divergence of the momentum conservation equation. Although this sets the reduced problem for the pressure, the definition of adequate boundary conditions remains an issue. Another way to construct the pressure reduced model is by duality of the momentum conservation equation with velocity supremizer test functions \cite{Kean2020}. 
\\
In this paper we afford the solution of the pressure equation by duality of the momentum conservation equation with gradients of the reduced pressure basis functions. This turns out to solving a Poisson equation for the pressure, but with the advantage of setting the right boundary conditions and thus ensuring an improved accuracy in its recovery. We shall call this method SM (Stabilisation-Motivated) ROM, as it  has been proposed in a framework of stabilised solution of incompressible flows \cite{IliescuJohn14}, see also \cite{Tallet2015a, Leblond2011} where the method is introduced from a minimum residual projection approach. We carry on the numerical analysis of the SM-ROM and prove the stability in specific norms. For velocity, in $L^\infty(L^2)$ and energy norms, with bounds that do not depend on the viscosity, and for pressure, in a semi-norm of the same asymptotic order as the $L^2$ norm with respect to the mesh size. The error estimates between the exact and the ROM-POD solutions are also obtained in the same norms. These estimates are optimal with respect to the projection error on the POD reduced spaces. To check the theoretical estimates, a numerical study is performed on two flow examples, the flow past a cylinder at $Re=200$ and the lid-driven cavity flow at $Re=9500$. The quality of the estimate is assessed by comparing the logarithmic slopes of the calculated versus the predicted errors with with respect to the velocity POD contribution ratio. It is shown on these examples that the error estimates allow a good approximation of the real errors slopes and thus a good prediction of the their rate of convergence. 
\\
The paper is organised as follows. In Section \ref{sec:VF} we introduce the discretization of the incompressible Navier-Stokes equations (NSE) that we will consider, while in Section \ref{sec:FOM} we describe its full order discretization by stable velocity-pressure finite elements. Section \ref{sec:POD-ROM} deals with the POD-ROM including the pressure recovery strategy that we are considering. Section \ref{sec:NA} is devoted to the numerical analysis, stability and error bounds. Numerical experiments on the flow past a cylinder and the lid-driven cavity flow are presented in Section \ref{sec:NumStud}. Finally conclusions are drawn in Section \ref{sec:Concl}
 }
\section{Time-dependent NSE: model problem and variational formulation}\label{sec:VF}

We introduce an Initial-Boundary Value Problem (IBVP) for the incompressible evolution NSE. 
For the sake of simplicity, we just impose the homogeneous Dirichlet BC on the whole boundary.

\medskip

Let $[0,T]$ be the time interval and $\Om$ a bounded polyhedral domain in $\mathbb{R}^{d}$, $d=2$ or $3$, with a Lipschitz-continuous boundary $\Ga=\partial\Om$. The transient NSE for an incompressible fluid are given by:

\medskip
\hspace{1cm} {\em Find $\uv:\Om\times (0,T)\longrightarrow\mathbb{R}^{d}$ and $p:\Om\times (0,T)\longrightarrow\mathbb{R}$ such that:}
\BEQ\label{eq:uNS}
\left \{
\begin{array}{rcll}
\partial_{t}\uv + (\uv \cdot\nabla)\uv - \nu \Delta\uv + \nabla p&=&\fv & \qmbx{in} \Om\times (0,T),\\
\div \uv &=&0 & \qmbx{in} \Om\times (0,T),\\
\uv &=& \bf{0} & \qmbx{on} \Ga\times (0,T),\\
\uv(\xv,0)&=&\uv_{0}(\xv) & \qmbx{in} \Om.
\end{array}
\right .
\EEQ
The unknowns are the velocity $\uv(\xv,t)$ and the pressure $p(\xv,t)$ of the incompressible fluid. The data are the source term $\fv(\xv,t)$, which represents a body force per mass unit (typically  the gravity), the kinematic viscosity $\nu$ of the fluid, which is a positive constant, and the initial velocity $\uv_{0}(\xv)$. 

\medskip

To define the weak formulation of problem (\ref{eq:uNS}), we need to introduce some useful notations for functional spaces \cite{Brezis11}. We consider the Sobolev spaces $H^{s}(\Om)$, $s\in \mathbb{R}$, $L^{p}(\Om)$ and $W^{m,p}(\Om)$, $m\in \mathbb{N}$, $1\leq p\leq\infty$. We shall use the following notation for vector-valued Sobolev spaces: ${\bf H}^{s}$, ${\bf L}^{p}$ and ${\bf W}^{m,p}$ respectively shall denote $[H^{s}(\Om)]^{d}$, $[L^{p}(\Om)]^{d}$ and $[W^{m,p}(\Om)]^{d}$ (similarly for tensor spaces of dimension $d\times d$). Also, the parabolic Bochner function spaces $L^{p}(0,T;X)$ and $L^{p}(0,T;{\bf X})$, where $X$ (${\bf X}$) stands for a scalar (vector-valued) Sobolev space, shall be denoted by $L^{p}(X)$ and $L^{p}({\bf X})$, respectively. In order to give a variational formulation of problem (\ref{eq:uNS}), let us consider the velocity space:
$$
\Xv={\bf H}_{0}^{1}=[H_{0}^{1}(\Om)]^{d}=\left\{\vv\in [H^{1}(\Om)]^{d} : \vv={\bf 0} \text{ on } \Ga\right\}. 
$$
This is a closed linear subspace of ${\bf H}^{1}$ and thus a Hilbert space endowed with the ${\bf H}^{1}$-norm. 
Thanks to Poincar\'e inequality, the ${\bf H}^{1}$-norm is equivalent on ${\bf H}_{0}^{1}$ to the norm $\nor{\vv}{{\bf H}_{0}^{1}}=\nor{\nabla\vv}{{\bf L}^{2}}$. Also, let us consider the pressure space:
$$
Q=L_{0}^{2}(\Om)=\left\{q\in L^{2}(\Om) : \int_{\Om}q\,d\xv=0\right\}. 
$$

\medskip

We shall thus consider the following variational formulation of \eqref{eq:uNS}:

\medskip
\hspace{1cm} {\em Given $\fv\in L^{2}({\bf L}^{2})$, find $\uv: (0,T)\longrightarrow\Xv$, $p: (0,T)\longrightarrow Q$ such that}
\BEQ\label{eq:fvNS} 
\left \{ 
\begin{array}{rcll}
\disp\frac{d}{dt}(\uv,\vv) + b(\uv,\uv,\vv) + \nu (\nabla\uv,\nabla\vv) - (p,\div\vv) &=& (\fv,\vv) & \forall\vv\in\Xv,\quad \text{in } \mathcal{D}'(0,T),\\
(\div \uv,q) &=&0 & \forall q \in Q, \quad \text{a.e. in } (0,T),\\
\uv(0)&=&\uv_{0},
\end{array} 
\right .
\EEQ
where $( \cdot,\cdot )$ stands for the $L^2$-inner product in $\Om$, and $\mathcal{D}'(0,T)$ is the space of distributions in $(0,T)$. The trilinear form
$b$ is given by: {\it for $\uv,\, \vv,\,\wv \in \Xv$}
\BEQ\label{eq:formb0}
b(\uv,\vv,\wv)={1 \over 2} \left[   ( \uv \cdot \nabla\, \vv, \wv) - ( \uv \cdot \nabla\,  \wv, \vv)\right].
\EEQ

\section{Finite element full order model}\label{sec:FOM}

In order to give a FE approximation of \eqref{eq:fvNS}, let $\{{\cal T}_{h}\}_{h>0}$ be a family of affine-equivalent, conforming (i.e., without hanging nodes) and regular triangulations of $\overline{\Om}$, 
formed by triangles or quadrilaterals ($d=2$), tetrahedra or hexahedra ($d=3$). For any mesh cell $K \in {\cal T}_{h}$,
its diameter will be denoted by $h_K$ and $h = \max_{K \in {\cal T}_{h}} h_K$. We consider $\Xv_{h}\subset\Xv$, $Q_{h}\subset Q$ being suitable FE spaces for velocity and pressure, respectively. The FE approximation of \eqref{eq:fvNS} can be written as follows:
\medskip

\hspace{1cm}{\em Find $(\uhv,p_h):(0,T)\longrightarrow \Xv_{h}\times Q_{h}$ such that}
\BEQ\label{eq:FEapprox} 
\left \{ 
\begin{array}{rcll}
\disp\frac{d}{dt}(\uhv,\vhv) + b(\uhv,\uhv,\vhv) + \nu (\nabla\uhv,\nabla\vhv) - (p_h,\div\vhv) &=& (\fv,\vhv) & \text{in } \mathcal{D}'(0,T),\\
(\div \uhv,q_h) &=&0 & \text{a.e. in } (0,T),\\
\uhv(0)&=&\uv_{0h},
\end{array} 
\right .
\EEQ
for any $(\vhv,q_h)\in \Xv_{h}\times Q_{h}$, and the initial condition $\uv_{0h}$ is some stable approximation to $\uv_{0}$ in $L^2$-norm belonging to $\Xv_{h}$.

\medskip

Actually, as a direct method we consider a Galerkin method with grad-div stabilization, since our aim is to get error bounds with constants independent on inverse powers of $\nu$. To describe this approach, we define hereafter the specific choice of FE spaces done both for the numerical analysis and practical computations in the present work. Given an integer $l\geq 2$ and a mesh cell $K \in {\cal T}_{h}$, denote by $\mathbb{R}^{l}(K)$ either $\mathbb{P}^{l}(K)$ (i.e., the space of Lagrange polynomials of degree $\leq l$, defined on $K$), if the grids are formed by triangles ($d=2$) or tetrahedra ($d=3$), or $\mathbb{Q}^{l}(K)$ (i.e., the space of Lagrange polynomials of degree $\leq l$ on each variable, defined on $K$), if the family of triangulations is formed by quadrilaterals ($d=2$) or hexahedra ($d=3$). We consider the mixed FE pair \review{known} as Taylor--Hood elements $(\Xv_{h}^l,Q_{h}^{l-1})$ \cite{BF,hood0} for the velocity-pressure:
\BEQ\label{eq:defVelSp}
\left \{ 
\begin{array}{lll}
&Y_{h}^{l}=V_{h}^{l}(\Om) = \{v_{h} \in C^{0}(\overline{\Om}) : v_{h}|_{K} \in \mathbb{R}^{l}(K),\,\forall K\in {\cal T}_{h}\},&\\ \\
&\Yhv^{l}=[Y_{h}^{l}]^{d} = \{\vhv \in [C^{0}(\overline{\Om})]^{d} : \vv_{h}|_{K}\in [\mathbb{R}^{l}(K)]^{d},\,\forall K\in {\cal T}_{h}\},&\\ \\
&\Xhv^l = \Yhv^{l}\cap {\bf H}_{0}^{1},\quad Q_h^{l-1}=Y_{h}^{l-1}\cap L_0^{2}.&
\end{array} 
\right .
\EEQ
Let us also consider the discrete space of divergence-free functions:
$$
\Vv_{h}^l=\left\{\vhv\in \Xhv^l : (\div\vhv,q_{h})=0\quad \forall q_{h}\in Q_{h}^{l-1}\right\}.
$$

\medskip

The considered grad-div FOM is given by: 
\medskip

\hspace{1cm}{\em Find $(\uhv,p_h):(0,T)\longrightarrow \Xhv^l\times Q_{h}^{l-1}$ such that}
\BEQ\label{eq:FOM} 
\left \{ 
\begin{array}{rcll}
\disp\frac{d}{dt}(\uhv,\vhv) + b(\uhv,\uhv,\vhv) + \nu (\nabla\uhv,\nabla\vhv) &-& (p_h,\div\vhv) &\\
+ \mu(\div\uhv,\div\vhv)&=&(\fv,\vhv) & \text{in } \mathcal{D}'(0,T),\\ \\
(\div \uhv,q_h) &=&0 & \text{a.e. in } (0,T),\\
\uhv(0)&=&\uv_{0h},
\end{array} 
\right .
\EEQ
for any $(\vhv,q_h)\in \Xhv^l\times Q_{h}^{l-1}$, where $\mu(\div\uhv,\div\vhv)$ denotes the grad-div stabilization term and $\mu$ is the positive grad-div stabilization parameter.

\medskip

To state the full space-time discretization of the unsteady grad-div FOM \eqref{eq:FOM}, consider a positive integer number $N$ and define $\Delta t = T/N$, $t_{n}=n\Delta t$, $n=0,1,\ldots,N$. We compute the approximations $\uhv^{n}$, $\ph^{n}$ to $\uv^{n}=\uv(\cdot,t_{n})$ and $p^{n}=p(\cdot,t_{n})$ by using, for simplicity of the analysis, a backward Euler scheme:
\begin{itemize}
\item {\bf Initialization.} Set: $\uhv^{0}=\uohv.$
\item {\bf Iteration.} For $n=0,1,\ldots,N-1$: {\em Given $\uhv^{n}\in\Xhv^l$, find $(\uhv^{n+1},\ph^{n+1})\in\Xhv^l\times Q_h^{l-1}$ such that:}
\BEQ\label{eq:discTempFOM} 
\left \{ 
\begin{array}{rcl}
\left(\disp\frac{\uhv^{n+1}-\uhv^{n}}{\Delta t},\vhv\right)+b(\uhv^{n+1},\uhv^{n+1},\vhv)+\nu(\nabla\uhv^{n+1},\nabla\vhv)& &\\
-(\ph^{n+1},\div\vhv) + \mu(\div\ModifMourad{\uhv^{n+1}},\div\vhv)&=& (\fv^{n+1}, \vhv),\\ \\
(\div\uhv^{n+1},\qh) &=& 0,
\end{array} 
\right .
\EEQ
for any $(\vhv,\qh) \in \Xhv^l\times Q_h^{l-1}$.
\end{itemize}

\medskip

It is well-known that considering the discrete divergence-free space $\Vv_{h}^l$ we can remove the pressure from \eqref{eq:discTempFOM} since $\uhv^{n+1}\in \Vv_{h}^l$ satisfies
\begin{eqnarray}\label{eq:discTempFOMdivfree}
&&\left(\disp\frac{\uhv^{n+1}-\uhv^{n}}{\Delta t},\vhv\right)+b(\uhv^{n+1},\uhv^{n+1},\vhv)+\nu(\nabla\uhv^{n+1},\nabla\vhv)
\nonumber\\
&&\quad+
\mu(\div\ModifMourad{\uhv^{n+1}},\div\vhv)=(\fv^{n+1}, \vhv), \quad \forall \vhv\in \Vv_{h}^l.
\end{eqnarray}
For this method the following bounds hold (see \cite{NS_grad_div,BoscoJulia19,BoscoJulia19Corr}):
\begin{eqnarray}\label{eq:cota_grad_div}
\|\bu^n-\bu_h^n\|_{{\bf L}^2}+ \left(\nu\sum_{j=1}^n\Delta t\|\nabla(\bu^j-\bu_h^j)\|_{{\bf L}^2}^2\right)^{1/2}\le C(\bu,p,l+1) (h^{l}+\Delta t), \quad 1\le n\le N,
\end{eqnarray}
and
\begin{eqnarray}\label{eq:cota_grad_div_pre}
\left(\sum_{j=1}^n\Delta t \|p^j-p_h^j\|_{{\bf L}^2}^2\right)^{1/2}\le C(\bu,p,l+1) h^{-1/2}(h^{l}+\Delta t), \quad 1\le n\le N,
\end{eqnarray}
where the constant $C(\bu,p,l+1)$ does not depend on inverse powers of $\nu$.

\section{Proper orthogonal decomposition reduced order model}\label{sec:POD-ROM}

We briefly describe the POD method, following \cite{KunischVolkwein01}, and apply it to the projection-based stabilized FOM \eqref{eq:discTempFOM}.

\medskip

Let us consider the ensembles of velocity snapshots $\chi^{v}=\text{span}\left\{\uhv^{1},\ldots,\uhv^{N}\right\}$ and pressure snapshots $\chi^{p}=\text{span}\left\{\ph^{1},\ldots,\ph^{N}\right\}$, given by the FE solutions to \eqref{eq:discTempFOM} at time $t_{n}$, $n=1,\ldots,N$. The POD method seeks low-dimensional bases $\left\{\boldsymbol{\varphi}_{1},\ldots,\boldsymbol{\varphi}_{r_v}\right\}$ and $\left\{\psi_{1},\ldots,\psi_{r_p}\right\}$ in real Hilbert spaces $\mathcal{H}_{v}$, $\mathcal{H}_{p}$ that optimally approximate the velocity and pressure snapshots with respect to the discrete $L^2(\mathcal{H}_{v})$, $L^2(\mathcal{H}_{p})$ norms, respectively (cf. \cite{KunischVolkwein01}). It can be shown that the following POD projection error formulas hold \cite{Holmes96, KunischVolkwein01}:
\BEQ\label{eq:PODerrVel}
\Delta t\sum_{n=1}^{N}\left\| \uhv^{n} - \sum_{i=1}^{r_v}\left(\uhv^{n},\boldsymbol{\varphi}_{i}\right)_{\mathcal{H}_{v}}\boldsymbol{\varphi}_{i} \right\|_{\mathcal{H}_{v}}^{2} = \sum_{i=r_{v}+1}^{M_{v}}\lambda_{i},
\EEQ
and
\BEQ\label{eq:PODerrPres}
\Delta t\sum_{n=1}^{N}\left\| \ph^{n} - \sum_{i=1}^{r_p}\left(\ph^{n},\psi_{i}\right)_{\mathcal{H}_{p}}\psi_{i} \right\|_{\mathcal{H}_{p}}^{2} = \sum_{i=r_{p}+1}^{M_{p}}\gamma_{i},
\EEQ
where $M_{v},M_{p}$ are the rank of $\chi^v$ and $\chi^p$, respectively, and $\lambda_{i},\gamma_{i}$ are the associated eigenvalues. Although $\mathcal{H}_{v},\mathcal{H}_{p}$ can be any real Hilbert spaces, in what follows we consider $\mathcal{H}_{v}={\bf L}^{2}$ and $\mathcal{H}_{p}=L^{2}$. Also, we are going to take the same number of velocity and pressure POD bases functions, i.e. $r_{v}=r_{p}=r$ in what follows.

\medskip

We respectively consider the following velocity and pressure spaces for the POD setting:
$$
\Xv_{r}=\text{span}\left\{\boldsymbol{\varphi}_{1},\ldots,\boldsymbol{\varphi}_{r}\right\}\subset \Xhv^l,
$$
and
$$
Q_{r}=\text{span}\left\{\psi_{1},\ldots,\psi_{r}\right\}\subset Q_h^{l-1}.
$$

\begin{remark}
Since the POD velocity modes are linear combinations of the snapshots, they satisfy the boundary conditions in \eqref{eq:uNS} and are solenoidal. Thus, the POD velocity modes belong to $\Vv_{h}^l$, which yields $\Xv_{r}\subset\Vv_{h}^l$.
\end{remark}

The standard Galerkin projection-based POD-ROM uses both Galerkin truncation and Galerkin projection.
The former yields an approximation of the velocity and pressure fields by a linear combination of the corresponding truncated POD basis:
\BEQ\label{eq:PODsolVel}
\uv(\xv,t)\approx \uv_{r}(\xv,t)=\sum_{i=1}^{r}a_{i}(t)\boldsymbol{\varphi}_{i}(\xv),
\EEQ
and
\BEQ\label{eq:PODsolVel}
p(\xv,t)\approx p_{r}(\xv,t)=\sum_{i=1}^{r}b_{i}(t)\psi_{i}(\xv),
\EEQ
where $\left\{a_{i}(t)\right\}_{i=1}^{r}$ and $\left\{b_{i}(t)\right\}_{i=1}^{r}$
are the sought time-varying coefficients representing the POD-Galerkin velocity and pressure trajectories. Note that $r << {\cal N}$, where ${\cal N}$ denotes the number of degrees of freedom (d.o.f.) in a full order simulation. Replacing the velocity-pressure FE pair $(\uhv,\ph)$ with $(\uv_{r},p_r)$ in the FE approximation \eqref{eq:discTempFOM} and projecting the resulted equations onto the POD product space $(\Xv_{r},Q_r)$ using the POD basis $\left(\left\{\boldsymbol{\varphi}_{i}\right\}_{i=1}^{r},\left\{\psi_{i}\right\}_{i=1}^{r}\right)$, the full space-time discretization of the grad-div ROM reads as: 
\begin{itemize}
\item {\bf Initialization.} Set: $\uv_{r}^{0}=\disp\sum_{i=1}^{r}(\uov,\boldsymbol{\varphi}_{i})\boldsymbol{\varphi}_{i}.$
\item {\bf Iteration.} For $n=0,1,\ldots,N-1$: {\em Given $\uv_{r}^{n}\in\Xv_{r}$, find $\uv_{r}^{n+1}\in\Xv_{r}$ such that:}
\begin{eqnarray}\label{eq:discTempROMvel}
&&\left(\disp\frac{\uv_{r}^{n+1}-\uv_{r}^{n}}{\Delta t},\boldsymbol{\varphi}\right)+b(\uv_{r}^{n+1},\uv_{r}^{n+1},\boldsymbol{\varphi})+\nu(\nabla\uv_{r}^{n+1},\nabla\boldsymbol{\varphi})
+\mu(\div\cb \uv_r^{n+1} \cn,\div\vv_r)\nonumber\\
&&\quad=(\fv^{n+1}, \boldsymbol{\varphi}),\quad \forall \boldsymbol{\varphi}\in  \Xv_{r}.
\end{eqnarray}
\end{itemize}
In \eqref{eq:discTempROMvel}, the pressure term vanishes due to the fact that $\uv_r^n$ belongs to the discrete divergence-free space $\Vv_{h}^l$. In order to recover the reduced order pressure, we propose to solve in a post-process phase the residual-based equation, as follows: \cb  Find $p_r^{n+1} \in Q_{r}$ such that 
\begin{eqnarray}\label{eq:discTempROMpres}
&&\sum_{K\in{\cal T}_h}\tau_K(\nabla p_r^{n+1},\nabla\psi)_K\\
&&\quad=-\sum_{K\in{\cal T}_h}\tau_K(\frac{\uv_r^{n+1}-\uv_r^{n}}{\Delta t}+ \tilde{\uv}_{r}^{n+1}\cdot\nabla\uv_r^{n+1}-\nu\,\Delta\uv_r^{n+1}-\fv^{n+1},\nabla\psi)_K,\nonumber
\end{eqnarray}
for all $\psi\in  Q_{r}$. Here, the $\tau_K$ are positive  coefficients of order $h_K^2$, that is, there exist two constants $C_1 >0$, $C_2 >0$, such that
\begin{equation}\label{estaus}
C_1 \, h_K^2 \le \tau_K \le C_2 \, h_K^2,\quad \forall K\in {\cal T}_h,\,\, \forall h>0.
\end{equation}
Also, the symbol $\tilde{\uv}_{r}$ stands for 
${\uv}_{r} $ when $d=2 $, and for some truncation of ${\uv}_{r}$ such that
\begin{equation} \label{hip:utrunc}
\|\tilde{\uv}_{r}\|_{L^\infty(K)} \le C_T \, h_K^{-1}, \quad \forall K\in {\cal T}_h,\,\, \forall h>0
\end{equation}
when $d=3$, for some constant $C_T>0$. We shall denote in this way all generic constants that will appear along the paper, unless needed to specify them by some index. Such a $\tilde{\uv}_{r}$ may be, for instance, defined as follows:
$$
\tilde{\uv}_{r} = \Pi_h(T_\alpha ({\uv}_{r})),\quad\mbox{with  } \alpha = h_K^{-1}
$$
where $\Pi_h$ denotes the standard Lagrange interpolation operator on space $\Vv_{h}^l$, $T_\alpha$ is the truncation at height $\alpha >0$, and $\alpha$ is the function defined on $\bar{\Omega}$ that assigns to any $x \in \bar{\Omega}$ the inverse of the average of the grid elements that contain $x$. It holds 
$$
T_\alpha(\uv)(x_i)= \left \{ \begin{array} {ccl}
\uv(x_i) & \mbox{if} & |\uv(x_i)| \le \alpha(x_i), \\
\alpha(x_i)\, sign(\uv(x_i)) & \mbox{if} & |\uv(x_i)| \ge \alpha(x_i).
\end{array}
\right .
$$
for any Lagrange interpolation node $x_i$. Then, as the grids are regular, $|T_\alpha(\uv)(x)|\le C h_K^{-1}$ for any $x \in \bar{\Omega}$, where $K\in \trh$ is the element containing $x$, for some constant $C>0$.

Due to inverse finite element error estimates and the local $L^2$ stability of the Lagrange interpolation (see \cite{chalew}, \review{Appendix),} 
$$\|\tilde{\uv}_{r}\|_{L^\infty(K)} \le C\, h_K^{-3/2}\|\tilde{\uv}_{r}\|_{L^2(K)} \le C\, h_K^{-3/2} \, \|T_\alpha({\uv}_{r})\|_{L^2(K)} \le C\, h_K^{-3/2} \, |K|^{1/2}\, \alpha\le  C\, h_K^{-1}.
$$
Note that the use of inverse finite element estimates  combined with the $L^\infty(L^2)$ bound for the $\uv_r$ yields
$$
\|{\uv}_{r}\|_{L^\infty(K)} \le C\, h_K^{-3/2}\|{\uv}_{r}\|_{L^2(K)} \le C\, h_K^{-3/2}.
$$
Then, the truncation at height $h_K^{-1}$, although somewhat restrictive, is rather mild, in particular it allows the blow-up of the $L^\infty(\Omega)$ norm of $\tilde{\uv}_{r}$ as $h \rightarrow 0$ . Anyhow, this technical hypothesis may be overcome when $d=2$ or a weaker norm of the pressure gradient is bounded when $d=3$.

The structure of method \eqref{eq:discTempROMpres} is inspired in that of stabilized methods for fluid flows. Usually a weighted residual (the difference between the l.h.s. and the r.h.s. of \eqref{eq:discTempROMpres}) is added to the Galerkin formulation to stabilize both convection dominance and pressure gradient \cite{metstab}. The scaling of the stabilizing coefficients $\tau_K$ as $h_K^2$ comes from dimensional analysis. Here, we isolate this stabilizing term to recover the pressure, once the velocity is known.

The presence of the factors $\tau_K$ in \eqref{eq:discTempROMpres} only is needed whenever the grids ${\cal T}_h$ are not uniformly regular. Otherwise all the $\tau_K$ may be taken equal to $h^2$ and \eqref{eq:discTempROMpres} may be replaced by 
\begin{eqnarray}\label{eq:discTempROMpres2}
(\nabla p_r^{n+1},\nabla\psi)_\Omega=-(\frac{\uv_r^{n+1}-\uv_r^{n}}{\Delta t}+ \tilde{\uv}_{r}^{n+1}\cdot\nabla\uv_r^{n+1}-\nu\,\Delta\uv_r^{n+1}-\fv^{n+1},\nabla\psi)_\Omega.
\end{eqnarray}

Note that equations \eqref{eq:discTempROMpres} and \eqref{eq:discTempROMpres2} contain a discretization of the natural b. c. for pressure, that is (assuming $\tilde{\uv}_{r}^{n+1}=\uv_{r}^{n+1}$)
$$
\nabla p_r^{n+1}\cdot {\bf n} =  \left (\frac{\uv_r^{n+1}-\uv_r^{n}}{\Delta t}+ \uv_r^{n+1}\cdot\nabla\uv_r^{n+1}-\nu\,\Delta\uv_r^{n+1}-\fv^{n+1}\right )\cdot {\bf n} \,\, \mbox{on   }\partial \Omega.
$$
Thus, no artificial boundary conditions on the pressure are being imposed.

 Hereafter, the grad-div ROM \eqref{eq:discTempROMvel} together with \eqref{eq:discTempROMpres} will be referred to as Stabilization Motivated-ROM (SM-ROM). \cb The SM-ROM \eqref{eq:discTempROMvel}-\eqref{eq:discTempROMpres}, without grad-div term, has been introduced and numerically investigated in \cite{IliescuJohn14}. 
 
 Let us remark that an alternative time discretization could be given by the semi-implicit Euler method, where the trilinear form in \eqref{eq:discTempFOM} is discretized by $b(\uhv^{n+1},\uhv^n,\vhv)$.
This semi-implicit time discretization is less costly from the computational point of view with respect to a fully implicit one, which yields a nonlinear algebraic system of equations to be solved. However, the numerical analysis will be performed in Section \label{sec:NA} for the more technical case of the fully implicit time discretization given by \eqref{eq:discTempFOM}. We shall present in Section \ref{sec:NumStud} numerical results obtained with the semi-implicit discretization.
\cn
\section{Analysis of the residual-based SM-ROM}\label{sec:NA}
In this section, we perform the numerical analysis of the proposed unsteady SM-ROM  \eqref{eq:discTempROMvel}-\eqref{eq:discTempROMpres}, dealing with stability and error estimates. 
\subsection{Technical background}\label{subsec:TB}
This section provides some technical results that are required for the numerical analysis. Throughout the paper, we shall denote by $C$ a positive constant that may vary from a line to another, but which is always independent of the viscosity $\nu$, the FE mesh size $h$, the FE velocity interpolation order $l$, the time step $\Delta t$, and the velocity, pressure eigenvalues $\lambda_i$, $\gamma_i$. 

We shall denote by $\|\cdot\|_{m,p}$ the norm in $W^{m,p}(\Omega)$ and by $\|\cdot\|_{m,p,K}$ the norm in $W^{m,p}(K)$ for some $K\in \trh$. Also, we shall denote by $\|\cdot\|_{m}$ the norm in $H^{m}(\Omega)$ and by $\|\cdot\|_{m,K}$ the norm in $H^{m}(K)$ for some $K\in \trh$. 

 We shall use some inverse estimates for finite element functions (cf. \cite{Ciarlet78}, Theorem 3.2.6):  there exists two constants $C_{I,1} >0$  and $C_{I,2} >0$ such that
\begin{equation}
\label{eq:invest}
\|v_h\|_{0,\infty,K} \le C_{I,1} \, h_K^{-\frac{d}{2}}\,\|v_h\|_{0,K}, \,\,\|\Delta v_h\|_{0,K} \le C_{I,2} \, h_K^{-1}\,\|\nabla v_h\|_{0,K}, \,\, \mbox{for all  } v_h \in V_h^l.
\end{equation}

We define the scalar product:
$$
(\cdot,\cdot)_{\tau}:{\bf L}^{2}(\Om)\times {\bf L}^{2}(\Om) \to \mathbb{R},\quad
(\vv,\wv)_{\tau} = \sum_{K\in{\cal T}_{h}}\tau_{K}(\vv,\wv)_{K},
$$
and its associated norm:
$$
\nor{\vv}{\tau}=(\vv,\vv)_{\tau}^{1/2}.
$$
We shall also use the following discrete Gronwall lemma:
\begin{lemma} \label{disgron}
Let $\alpha_n$, $\beta_n$, $\gamma_n$, $n=1,2,...$ be non-negative real numbers satisfying
\BEQ \label{estaux1}
(1-\sigma \,\Delta t) \,\alpha_{n+1} + \beta_{n+1} \le (1+\tau \,\Delta t) \,\alpha_n + \gamma_{n+1} 
\EEQ
for some $\sigma \ge 0$, $\tau \ge 0$. Assume that $\sigma \, \Delta t \le 1-\delta$ for some $\delta >0$. Then it holds
\BEQ \label{estgr1}
\alpha_n \le e^{\rho \,t_n}\,\alpha_0 + \frac{1}{\delta} \,\sum_{l=1}^n e^{\rho\, (t_n-t_l)}\, \gamma_l;
\EEQ
\BEQ \label{estgr2}
\sum_{l=1}^n\beta_l \le \left (1+\frac{\tau}{\sigma} +(\sigma+\tau)\,e^{\rho\, t_{n-1}}\,t_{n-1} \right )\,\alpha_0 +\frac{1}{\delta} \, \left ( 1+ (\sigma+\tau)\,e^{\rho\, t_{n-1}}\,t_{n-1} \right )\sum_{l=1}^n\,\gamma_l,
\EEQ
with $\disp \rho=\frac{\sigma+\tau}{\delta}$.
\end{lemma} \label{le:bounds}
\begin{proof} Dividing by $1-\sigma \,\Delta t$ in  \eqref{estaux1}, 
\begin{eqnarray}
\alpha_{n+1} +\beta_{n+1} &\le& \alpha_{n+1} +\frac{1}{1-\sigma\Delta t} \beta_{n+1} \le (1+\frac{(\sigma+\tau) \Delta t}{1-\sigma \Delta t})\alpha_n +\frac{1}{1-\sigma \Delta t} \,\gamma_{n+1} \nonumber \\
&\le& (1+\rho \Delta t)\alpha_n +\frac{1}{\delta} \,\gamma_{n+1}.
\end{eqnarray}
Recursively using this inequality,
\begin{eqnarray}
\alpha_{n+1} +\beta_{n+1} &\le& (1+\rho \Delta t) ^2\,\alpha_{n-1} +(1+\rho \Delta t) \,\, \frac{1}{\delta} \,\gamma_n +\frac{1}{\delta} \,\gamma_{n+1}\le \cdots \nonumber \\
&\le& (1+\rho \Delta t) ^{n+1}\,\alpha_0+\frac{1}{\delta} \,\sum_{l=1}^{n+1}(1+\rho \Delta t) ^{n+1-l} \,\gamma_l\nonumber \\
&\le& e^{\rho t_{n+1}}\,\alpha_0 +\frac{1}{\delta} \,\sum_{l=1}^{n+1}e ^{\rho (t_{n+1}-t_l)} \,\gamma_l,
\end{eqnarray}
hence \eqref{estgr1}. Summing now the inequalities \eqref{estaux1} yields
$$
(1-\sigma \,\Delta t)\, \sum_{l=1}^{n+1}\alpha_i+ \sum_{l=1}^{n+1}\beta_l \le (1+\tau \,\Delta t)\,\sum_{l=0}^n\alpha_i +\frac{1}{\delta} \,\sum_{l=1}^{n+1}\gamma_l .
$$
Then,
\begin{eqnarray}
\sum_{l=1}^{n+1}\beta_l &\le& (1+\tau \,\Delta t)\,\alpha_0+(\sigma+\tau)\Delta t\,\sum_{l=1}^{n}\alpha_l+\frac{1}{\delta} \,\sum_{l=1}^{n+1}\gamma_l \nonumber \\
&\le& (1+\tau \,\Delta t)\,\alpha_0+(\sigma+\tau)\, \Delta t\,\sum_{l=1}^{n}\left [e^{\rho t_l}\,\alpha_0 +\frac{1}{\delta} \,\sum_{r=1}^l e ^{\rho (t_l-t_r)} \,\gamma_l \right ] +\frac{1}{\delta} \, \sum_{l=1}^{n+1}\gamma_l  \nonumber \\
&\le& \left ( 1+\frac{\tau}{\sigma}+ (\sigma+\tau)\,e^{\rho t_n}\,t_n\right ) \, \alpha_0+\frac{1}{\delta}\, \left ( (\sigma+\tau)\,e^{\rho\,t_n}\,t_n+1 \right )\, \sum_{l=1}^{n+1}\gamma_l,
\end{eqnarray}
hence \eqref{estgr2}.
\end{proof}
Note that the constant $\rho$ plays the role of an amplification factor for the exponential dependence on the length of the time interval of the bounds appearing in \eqref{estgr1} and \eqref{estgr2}. It tends to 1 as $\Delta t$ tends to zero.
\subsection{Existence and stability results for SM-ROM}\label{subsec:STAB}

We have the following existence and unconditional stability result for the SM-ROM,\eqref{eq:discTempROMvel}-\eqref{eq:discTempROMpres}: \cb

\begin{theorem} Assume that $\Delta t \le 1-\delta$ for some $\delta >0$, and that $\fv \in L^2(0,T;{\bf L}^2(\Omega))$. Assume also that the data $\fn$ in problems \eqref{eq:discTempROMvel} and \eqref{eq:discTempROMpres} satisfy
$$
\|\uv^0_r\|_0 \le \|\uv_0\|_0,\quad \Delta t\, \sum_{i=1}^n \|\fn\|_0^2 \le \|\fv\|^2_{L^2(0,t_n;{\bf L}^2(\Omega))},\quad \mbox{for   } n=1,\cdots, N.
$$
Then the following statements hold.
\begin{enumerate}
\item Equations \eqref{eq:discTempROMvel} admits at least a solution $\uv_{r}^{k+1}$, and \eqref{eq:discTempROMpres} admits a unique solution $p_r^{n+1}$. 
\item Any solution $\uv_{r}^{k+1}$ of equation  \eqref{eq:discTempROMvel}  satisfies the following estimates
\begin{eqnarray} 
\|\uv_{r}^k\|_0^2   \le e^{\rho\, t_k}\, \|\uv_0\|_0^2 + \frac{1}{\delta}\, e^{\rho\, t_k}\, \|{\bf f}\|_{L^2(0,t_k;{\bf L}^2(\Omega)^2)}^2 , \quad k=1,\cdots, N;  \label{estvel1}&&
\end{eqnarray}
\begin{eqnarray} \label{estvel2}
\Delta t \, \sum_{n=1}^k \left (\nu\, \|\nabla \uv_{r}^{n}\|_0^2 + \mu \, \|\nabla\cdot  \uv_{r}^{n}\|_0^2 \right )  \le  \frac{1}{2}\, \left (1+t_k\,e^{\rho\, t_k}\right ) \|\uv_0\|_0^2  &&\\
+ \frac{1}{2\delta}\,  \left (1 + t_k\,e^{\rho\, t_k}\right )\, \|{\bf f}\|_{L^2(0,t_k;{\bf L}^2(\Omega))}^2 , \quad k=1,\cdots, N\nonumber &&
\end{eqnarray}
with $\disp \rho=\frac{1}{\delta}$.

\item Assume in addition that $h\le C_S\, \Delta t$ for some constant $C_S>0$, $\Delta t \le \disp \frac{1-\delta}{s+1/5}$ with $s=5\,C_2\, C_S^2$ for some $0<\rho<1$, and that in problem \eqref{eq:discTempROMpres} , 
\begin{itemize}
\item[3a)]  When $d=2$, it holds 
\begin{equation}\label{cond2}
C_2\le \disp \sqrt{\frac{\nu}{10}}\, \frac{1}{A\,C_{I,1}}, \,\mbox{ with  }
A=e^{\rho\, T}\, \|\uv_0\|_0^2 + \frac{1}{\delta}\, e^{\rho\, T}\, \|{\bf f}\|_{L^2(0,T;{\bf L}^2(\Omega)^2)}^2,
\end{equation}
or

\item[3b)]  When $d=3$, it holds $\|\tilde{\uv}_r^n\|_{0,\infty,K}\le C_T\, h_K^{-1}$ with 
\begin{equation}\label{cond3}
C_T \le \nu \,\sqrt{C_{I,2}},\quad C_2 \le \disp \frac{1}{10}\, \frac{1}{\nu\,C_{I,2}}.
\end{equation}
\end{itemize}
Then the solution $p_r^1,\cdots,p_r^N$ of equation \eqref{eq:discTempROMpres} satisfies the following estimates:
\begin{eqnarray}\label{estpress1}
\Delta t\,\sum_{n=1}^{N}\|\nabla \pn\|_\tau^2 &\le& \left [2+\left (2s+1/5\right ) \,e^{\rho t}\,T \right ] \|\uv_0\|_0^2  \\
&&+\frac{5}{\delta}\, \left [1+\left (2s+1/5\right)\, e^{\rho t}\,T \right ] \,(1+C_2\, h^2)\,\|\fv\|_{L^2(0,T;{\bf L}^2(\Omega))} .\nonumber
\end{eqnarray}
\end{enumerate}
\end{theorem}
\begin{proof} {\bf 1.-} Once $\uvn$ and $\uvnu$ are known, problem \eqref{eq:discTempROMpres} is a standard elliptic problem with a form $a_{pres}(p,q)=(\nabla p,\nabla q)_\Omega$ which is coercive on sub-spaces of $H^1(\Omega)$ with vanishing mean in $\Omega$. Then, problem \eqref{eq:discTempROMpres} admits a unique solution.

The existence of solution of equation \eqref{eq:discTempROMvel} follows from a standard finite-dimensional compactness argument, that we omit for brevity.

\vspace*{.3cm}\noindent
{\bf 2.-} Estimate \eqref{estvel2} follows from a standard argument for finite element approximation of incompressible flow problems using the discrete Gronwall Lemma \ref{disgron} that can be found, for instance, in \cite{chalew}, Chapter 10. We omit it for brevity

\vspace*{.3cm}\noindent
{\bf 3.-} Let us take $\vv_r = \uvnu$ in problem \eqref{eq:discTempROMvel}. Using the identities
$$
(\uvnu-\uvn,\uvnu)_\Om = \frac{1}{2}\,\|\uvnu\|_0^2 -\frac{1}{2}\,\|\uvn\|_0^2 + \frac{1}{2}\,\|\uvnu-\uvn\|_0^2,
$$
and
$$
b(\uvn;\uvnu,\uvnu)=0
$$
it holds
\begin{eqnarray} 
\frac{1}{2}\,\|\uvnu\|_0^2 -\frac{1}{2}\,\|\uvn\|_0^2 + \frac{1}{2}\,\|\uvnu-\uvn\|_0^2 + \nu\, \Delta t  \, \|\nabla \uvnu\|_0^2+ 2\, \mu\, \Delta t  \, \|\nabla \cdot \uvnu\|_0^2 &&\nonumber \\
 = \Delta t \, (\fnu,\uvnu)_\Omega\le \frac{\Delta t}{2}\, \| \fnu\|_0^2 +  \frac{\Delta t}{2}\, \| \uvnu\|_0^2. &&\label{estvelbasica}
\end{eqnarray}

Taking  $\psi=\pnu$ and multiplying by $\Delta t$ in problem \eqref{eq:discTempROMpres}, and summing up the first identity in \eqref{estvelbasica} yields
\begin{eqnarray}
&&\frac{1}{2}\,\|\uvnu\|_0^2 -\frac{1}{2}\,\|\uvn\|_0^2 + \nu\, \Delta t  \, \|\nabla \uvnu\|_0^2+\Delta t \, \|\nabla\pnu\|_\tau^2 \\
&& \le -\sum_{K\in \trh} \tau_K\,(\uvnu-\uvn+\Delta t\, (\utnu\cdot\nabla \uvnu -\nu\, \Delta \uvnu-\fnu),\nabla \pnu)_K  \nonumber \\
&&+ \Delta t \, (\fnu,\uvnu)_\Omega\nonumber
\end{eqnarray}
We bound each term on the r.h.s. of this expression. Observe that
$$
\sum_{K\in \trh} \tau_K\,(a,b)_K \le \|a\|_\tau \, \|b\|_\tau\le \frac{1}{2}\beta^{-1}\,  \|a\|_\tau^2+  \frac{1}{2}\beta\,\|b\|_\tau^2\quad
$$
for any functions $a,\, b \in L^2(\Omega)$ or $a,\, b \in {\bf L}^2(\Omega)$ and any $\beta >0$. Using $\tau_K \le C_2\, h_K^2$, $h_K \le C_S \,\Delta t$ and estimate \eqref{estvel1} it follows
\begin{eqnarray}
\sum_{K\in \trh} \tau_K\,(\uvnu,\nabla \pnu)_K \le \frac{1}{2\varepsilon\Delta t}\,\|\uvnu\|_\tau^2 + \frac{\varepsilon}{2}\,\Delta t\, \|\nabla \pnu\|_\tau^2 \nonumber  &&\\
\le    \frac{1}{2\varepsilon}\, \Delta t\, C_2\, C_S^2 \, \|\uvnu\|_0^2  + \frac{\varepsilon}{2}\,\Delta t\,\|\nabla \pnu\|_\tau^2; \label{estp1} &&
\end{eqnarray}
where $\varepsilon>0$ is a constant to be determined later. Similarly,
\begin{eqnarray}\label{estp2}
\sum_{n=0}^{N-1}\sum_{K\in \trh} \tau_K\,(\uvn,\nabla \pnu)_K  \le  \frac{1}{2\varepsilon}\, C_2\, C_S^2 \, \|\uvn\|_0^2   + \frac{\varepsilon}{2}\,\Delta t\,\sum_{n=1}^N\|\nabla \pn\|_\tau^2; &&
\end{eqnarray}
 Further,
\begin{eqnarray}
-\Delta t\,\sum_{K\in \trh} \tau_K\, (\nu\, \Delta\uvnu,\nabla \pnu)_K \le \frac{\nu^2}{2\varepsilon}\,\Delta t\,\|\Delta \uvnu\|_\tau^2 + \frac{\varepsilon}{2}\,\Delta t\, \|\nabla \pnu\|_\tau^2 \nonumber  &&\\
\le \frac{\nu^2}{2\varepsilon}\,C_{I,2}\, \max_{K\in \trh} \tau_K\,h_K^{-2} \,\Delta t\, \|\nabla \uvnu\|_0^2 + \frac{\varepsilon}{2}\,\Delta t\, \|\nabla \pnu\|_\tau^2 \nonumber &&\\
\le \frac{\nu^2}{2\varepsilon}\,C_{I,2}\, C_2 \,\Delta t\, \|\nabla \uvnu\|_0^2 + \frac{\varepsilon}{2}\,\Delta t\, \|\nabla \pnu\|_\tau^2.\label{estp3}
\end{eqnarray}
The convection-pressure gradient interaction term will be specifically bounded in each of the situations 3a) and 3b), as follows.

\vspace*{.3cm}\noindent
{\bf 3a)}  In this case $d=2$ and $\utnu=\uvnu$. It holds
\begin{eqnarray}
&\disp \Delta t \sum_{K\in \trh} \tau_K\,(\utnu\cdot\nabla \uvnu ,\nabla \pnu)_K \le \Delta t \sum_{K\in \trh} \tau_K\,\|\uvnu\|_{0,\infty,K}\, \|\nabla \uvnu\|_{0,K} \,\|\nabla \pnu\|_{0,K}  & \nonumber \\
&\disp \le  C_{I,1}\, \,\Delta t \sum_{K\in \trh} \tau_K\, h_K^{-d/2}\, \|\uvnu\|_{0,K}\,\,\|\nabla \uvnu\|_{0,K} \,\|\nabla \pnu\|_{0,K}& \nonumber \\
&\disp \le A\, C_{I,1}\, \,\Delta t \sum_{K\in \trh} \tau_K\, h_K^{-d/2}\,\|\nabla \uvnu\|_{0,K} \,\|\nabla \pnu\|_{0,K} ,&\nonumber\\
&\disp \le \frac{1}{2\varepsilon}\, A^2 \, C_{I,1}^2  \,\Delta t \, \sum_{K\in \trh} \tau_K\, h_K^{-d}\,\|\nabla \uvnu\|_{0,K}^2+ \frac{\varepsilon}{2}\,\Delta t\,\|\nabla \pnu\|_\tau^2 &\nonumber\\
&\disp \le \frac{1}{2\varepsilon} \, A^2 \, C_{I,1}^2\,C_2\,\Delta t\,\|\nabla \uvnu\|_{0,K}^2 + \frac{\varepsilon}{2}\,\Delta t\, \|\nabla \pnu\|_\tau^2. \label{estp4a}&
\end{eqnarray}
where $A$ is defined in \eqref{cond2}.

\vspace*{.3cm}\noindent
{\bf 3b)} In this case $d=3$ and $\|\utnu \|_{0,\infty,K}\le C_T \, h_K^{-1}$. It holds
\begin{eqnarray}
&\disp \Delta t \sum_{K\in \trh} \tau_K\,(\utnu\cdot\nabla \uvnu ,\nabla \pnu)_K \le \Delta t \sum_{K\in \trh} \tau_K\,\|\utnu\|_{0,\infty,K}\, \|\nabla \uvnu\|_{0,K} \,\|\nabla \pnu\|_{0,K}  & \nonumber \\
&\disp \le  C_T \,\Delta t \sum_{K\in \trh} \tau_K\, h_K^{-1}\,\|\nabla \uvnu\|_{0,K} \,\|\nabla \pnu\|_{0,K}& \nonumber\\
&\disp \le \frac{1}{2\varepsilon} \,C_T^2 \,\Delta t \, \sum_{K\in \trh} \tau_K\, h_K^{-2}\,\|\nabla \uvnu\|_{0,K}^2+ \frac{\varepsilon}{2}\,\Delta t\,\|\nabla \pnu\|_\tau^2 &\nonumber\\
&\disp \le \frac{1}{2\varepsilon} \,C_T^2 \,C_2\, \|\nabla \uvnu\|_{0}^2 + \frac{\varepsilon}{2}\,\Delta t\,\|\nabla \pnu\|_\tau^2. &\label{estp4b}
\end{eqnarray}
Further,
\begin{eqnarray}
\disp \Delta t \sum_{K\in \trh} \tau_K\,(\fnu ,\nabla \pnu)_K \le 
\disp \frac{1}{2\varepsilon} \,\Delta t\,\|\fnu\|_\tau^2 + \frac{\varepsilon}{2}\,\Delta t\,\|\nabla \pnu\|_\tau^2 \label{estp5};
\end{eqnarray}
and
\begin{eqnarray}
\disp \Delta t\,(\fnu ,\uvnu) \le 
\disp \frac{1}{2\varepsilon} \,\Delta t\,\|\fnu\|_0^2 + \frac{\varepsilon}{2}\,\Delta t\,\|\nabla \uvnu\|_0^2 \label{estp6};
\end{eqnarray}
Taking now $\varepsilon=1/5$ and summing up \eqref{estp1}, \eqref{estp2}, \eqref{estp3}, \eqref{estp5} and \eqref{estp6} with \eqref{estp4a} (case 3a)), or \eqref{estp4b} (case 3b)) yields
\begin{eqnarray} 
(1-(s+1/5)\Delta t)\,\|\uvnu\|_0^2   + \mu\, \Delta t  \, \|\nabla \uvnu\|_0^2 + \Delta t  \, \|\nabla \pnu\|_\tau^2\nonumber &&\\\le (1+s\Delta t)\, \|\uvn\|_0^2+ 5\,\Delta t\, \left (\| \fnu\|_0^2+ \| \fnu\|_\tau^2\right), &&\nonumber\\\le (1+s\Delta t)\, \|\uvn\|_0^2+ 5\,\Delta t\, (1+C_2\,h^2)\,\| \fnu\|_0^2 \label{estpresf} &&
\end{eqnarray}
where
$$\mu = 2\,\left ( \nu - \frac{1}{2\varepsilon} \, A^2\, C_{I,1}^2\,C_2 -  \frac{\nu^2}{2\varepsilon} \,C_2\, C_{I,2}\right ) \,\quad \mbox{when  } d=2,
$$
or
$$\mu = 2\,\left ( \nu - \frac{1}{2\varepsilon} \, C_2\, C_T^2 -  \frac{\nu^2}{2\varepsilon} \, C_2\, C_{I,2}\right ) \,\quad \mbox{when  } d=3.
$$
In any case, in view of conditions \eqref{cond2} and \eqref{cond3}, it holds 
$\mu \ge \nu$. Applying Lemma \ref{disgron} to \eqref{estpresf} with $\sigma = s+1/5$, $\tau =s$, $\alpha_n= \mu\, \Delta t  \, \|\uvnu\|_0^2 $, $\beta_n=\Delta t  \, \|\nabla \pnu\|_\tau^2$ and $\gamma_n=5\,\Delta t\, (1+h^2)\,\| \fnu\|_0^2$ yields estimate \eqref{estpress1} .
\end{proof}
\begin{remark} Estimates for velocity \eqref{estvel1}, \eqref{estvel2} and for pressure \eqref{estpress1} do not degenerate as $\nu \to 0$. This is obtained using the enhanced regularity $\fv \in L^2(0,T;{\bf L}^2(\Omega))$. Also, due to the conditions satisfied by constant $C_2$ in either \eqref{cond2} or \eqref{cond3}, \eqref{estpress1} provides either an estimate for the norm
$$
\nu\,\Delta t\, \sum_{i=1}^{N-1} \| p_r^{n+1}\|^2_h\,\,\,\mbox{when  } d=2, 
\,\,\mbox{or  } 
\Delta t\, \sum_{i=1}^{N-1} \| p_r^{n+1}\|^2_h,\,\mbox{when  } d=3,\,\,\mbox{with  } 
$$
$$
\|p\|_h = \left (\sum_{i=1}^{N-1} h_K^2\,\|\nabla p\|^2_{0,K}\right )^{1/2}.
$$
 The $\|\cdot\|_h$ norm plays a role in the discrete inf-sup condition for non-stable pairs of finite element velocity-pressure spaces (such as $\Xhv^l$ and $Q_{h}^l$, for instance, with equal interpolation in velocity and pressure, cf. \cite{Cha_maca}):
$$
\|p_h\|_0 \le \sup_{\vh \in \Xhv^l} \frac{(\div\vh,p_h)}{\|\vh\|_{1,\Omega}} + C\, \|p_h\|_h \,\mbox{for any  } \,\,p_h \in Q_h^l.
$$
Also, a straightforward argument using inverse inequalities proves that 
$$\|p_h\|_h \le C\, \| p_h\|_0\,\,\mbox{for any  } \,\,p_h \in Q_h^{l-1}.
$$
\end{remark}

\cn

\subsection{Error estimates for SM-ROM}\label{subsec:EE}
We are now in position to state the following error estimate results for the SM-ROM defined by \eqref{eq:discTempROMvel}-\eqref{eq:discTempROMpres}:
\begin{Theorem}[Velocity error estimate]
Let $\bu$ be the velocity in the NSE \eqref{eq:uNS}, let $\bu_r$ be the grad-div ROM
velocity defined in \eqref{eq:discTempROMvel}, and assume that the solution $(\bu,p)$ of \eqref{eq:uNS} is regular enough. Then, the
following bound holds:
\begin{equation}\label{eq:cota_finalSUPv}
\Delta t\sum_{n=0}^{N-1} \|\bu_r^{n+1}-\bu^{n+1}\|_{0}^2\le C\left( \sum_{i=r+1}^{M_v} \lambda_i \left(1+\nor{\nabla \boldsymbol{\varphi}_{i}}{0}^2\right)+h^{2l}+\Delta t^2\right).
\end{equation}
\end{Theorem}

The proof of this result follows along the same lines as the proof of Theorem 5.3 in \cite{NovoRubinoSINUM21}, we omit it for brevity.

\begin{Theorem}[Pressure error estimate]
Let $p$ be the pressure in the NSE \eqref{eq:uNS}, let $p_r$ be the SM-ROM
pressure defined in \eqref{eq:discTempROMpres}, and assume that the solution $(\bu,p)$ of \eqref{eq:uNS} is regular enough. Then, the
following bound holds:
\begin{equation}\label{eq:cota_finalSUPp}
\Delta t\sum_{n=0}^{N-1} \|\nabla (p_r^{n+1}-p^{n+1})\|_{\tau}^2\le C E^{*},
\end{equation}
where: 
\begin{eqnarray}\label{eq:Estar}
E^{*}&=&{h^2}+(h^{2l}+\Delta t^2)(1+\nu^{-1}+\nor{S_r^v}{2}+h^{-1})\nonumber \\
&&+\sum_{i=r+1}^{M_v}\lambda_i\left[1+\left(1+\nor{S_r^v}{2}\right)\nor{\nabla\boldsymbol{\varphi}_i}{0}^2\right]+h^2\sum_{i=r+1}^{M_p}\gamma_i\nor{\nabla\psi_i}{0}^2,
\end{eqnarray}
{and the penalty term $h^2$ in \eqref{eq:Estar} actually vanishes when $d=2$.
}
\end{Theorem}

{\bf Proof.} We start deriving the pressure error bound by splitting the error for the velocity and the pressure into two terms:
\BEQ\label{eq:splitErrVel}
\uv^{n+1}-\uv_{r}^{n+1}=(\uv^{n+1}-P_{r}^v \uv^{n+1})-(\uv_{r}^{n+1}-P_{r}^v \uv^{n+1})=\boldsymbol{\eta}^{n+1}-\boldsymbol{\phi}_{r}^{n+1},
\EEQ
\BEQ\label{eq:splitErrPres}
p^{n+1}-p_{r}^{n+1}=(p^{n+1}-P_r^p p^{n+1})-(p_{r}^{n+1}-P_r^p p^{n+1})=\rho^{n+1}-s_{r}^{n+1}.
\EEQ
In \eqref{eq:splitErrVel}, the first term, $\boldsymbol{\eta}^{n+1}=\uv^{n+1}-P_{r}^v \uv^{n+1}$, represents the difference between $\uv^{n+1}$ and its $L^2$-orthogonal projection on $\Xv_r$. The second term, $\boldsymbol{\phi}_{r}^{n+1}=\uv_{r}^{n+1}-P_{r}^v \uv^{n+1}$, is the remainder. Similarly, in \eqref{eq:splitErrPres}, the first term, $\rho^{n+1}=p^{n+1}-P_r^p p^{n+1}$, represents the difference between $p^{n+1}$ and its $L^2$-orthogonal projection on $Q_r$. The second term, $s_{r}^{n+1}=p_{r}^{n+1}-P_r^p p^{n+1}$, is the remainder. 

\medskip

Next, we construct the error equation. We first consider the solution $(\bu,p)$ of \eqref{eq:uNS} at $t=t_{n+1}$ and:
\BEQ\label{eq:ContEq}
\sum_{K\in{\cal T}_h}\tau_K(\partial_t \uv^{n+1}+ \uv^{n+1}\cdot\nabla\uv^{n+1}-\nu\,\Delta\uv^{n+1}+\nabla p^{n+1}-\fv^{n+1},\nabla s_{r}^{n+1})_K=0,
\EEQ
then subtract \eqref{eq:discTempROMpres} with $\psi=s_{r}^{n+1}$ from it. By adding and subtracting the difference quotient term $\disp\frac{1}{\Delta t}(\uv^{n+1}-\uv^{n})$, and applying the decompositions \eqref{eq:splitErrVel}-\eqref{eq:splitErrPres}, we get:
\BEQ\label{eq:ErrEq0}
\begin{array}{lll}
&&\disp\sum_{K\in{\cal T}_h}\tau_K\left(\partial_{t}\uv^{n+1}-\disp\frac{\uv^{n+1}-\uv^{n}}{\Delta t}+\disp\frac{1}{\Delta t}(\boldsymbol{\eta}^{n+1}-\boldsymbol{\phi}_{r}^{n+1})
-\disp\frac{1}{\Delta t}(\boldsymbol{\eta}^{n}-\boldsymbol{\phi}_{r}^{n}),\nabla s_r^{n+1}\right)_{K}\\ \\
&+&\disp\sum_{K\in{\cal T}_h}\tau_K\left( (\uv^{n+1}\cdot\nabla\uv^{n+1}-\tilde{\uv}_r^{n+1}\cdot\nabla\uv_r^{n+1})
-\nu\Delta(\boldsymbol{\eta}^{n+1}-\boldsymbol{\phi}_{r}^{n+1}),\nabla s_r^{n+1}\right)_{K}\\ \\
&+&\disp\sum_{K\in{\cal T}_h}\tau_K\left(\nabla(\rho^{n+1}-s_r^{n+1}),\nabla s_r^{n+1}\right)_{K}=0.
\end{array}
\EEQ
Letting $\boldsymbol{c}^{n}=\partial_{t}\uv^{n+1}-\disp\frac{\uv^{n+1}-\uv^{n}}{\Delta t}$, from \eqref{eq:ErrEq0} we obtain:
\BEQ\label{eq:ErrEq1}
\begin{array}{lll}
&&\Delta t\disp\sum_{K\in{\cal T}_h}\tau_K\nor{\nabla s_r^{n+1}}{0,K}^2 = \disp\sum_{K\in{\cal T}_h}\tau_K
\left(\Delta t \boldsymbol{c}^{n} + (\boldsymbol{\eta}^{n+1}-\boldsymbol{\eta}^{n}) - (\boldsymbol{\phi}_{r}^{n+1}-\boldsymbol{\phi}_{r}^{n}),\nabla s_r^{n+1}\right)_K\\ \\
&+&\Delta t\disp\sum_{K\in{\cal T}_h}\tau_K\left( (\uv^{n+1}\cdot\nabla\uv^{n+1}-\tilde{\uv}_r^{n+1}\cdot\nabla\uv_r^{n+1})
-\nu\Delta(\boldsymbol{\eta}^{n+1}-\boldsymbol{\phi}_{r}^{n+1}),\nabla s_r^{n+1}\right)_{K}\\ \\
&+&\Delta t\disp\sum_{K\in{\cal T}_h}\tau_K\left(\nabla\rho^{n+1},\nabla s_r^{n+1}\right)_{K}=I+II+III+IV+V+VI.
\end{array}
\EEQ
We estimate the terms on the r.h.s. of \eqref{eq:ErrEq1} one by one. By Cauchy-Schwarz and Young' s inequalities, we bound the first term on the r.h.s. of \eqref{eq:ErrEq1}:
\BEQ\label{eq:cotaI}
I=\Delta t\sum_{K\in{\cal T}_h}\tau_K (\boldsymbol{c}^{n},\nabla s_r^{n+1})_K\leq \Delta t \sum_{K\in{\cal T}_h} \frac{\tau_K}{2\varepsilon}\nor{\boldsymbol{c}^{n}}{0,K}^2 + \Delta t \sum_{K\in{\cal T}_h}\tau_K \frac{\varepsilon}{2}\nor{\nabla s_r^{n+1}}{0,K}^2,
\EEQ
for some small positive constant $\varepsilon$ (to be determined later).

\medskip

For the second term on the r.h.s. of \eqref{eq:ErrEq1}, assuming $h\sim \Delta t$, we have:
\begin{eqnarray}\label{eq:cotaII}
II&=&-\sum_{K\in{\cal T}_h}\tau_K (\boldsymbol{\eta}^{n+1}-\boldsymbol{\eta}^{n},\nabla s_r^{n+1})_K\nonumber \\
&\leq&  \sum_{K\in{\cal T}_h} \frac{\tau_K}{2\varepsilon\Delta t}\nor{\boldsymbol{\eta}^{n+1}-\boldsymbol{\eta}^{n}}{0,K}^2 + \Delta t \sum_{K\in{\cal T}_h}\tau_K \frac{\varepsilon}{2}\nor{\nabla s_r^{n+1}}{0,K}^2\nonumber \\
&\leq& \Delta t \sum_{K\in{\cal T}_h} \frac{C}{\varepsilon}\nor{\boldsymbol{\eta}^{n+1}-\boldsymbol{\eta}^{n}}{0,K}^2 + \Delta t \sum_{K\in{\cal T}_h}\tau_K \frac{\varepsilon}{2}\nor{\nabla s_r^{n+1}}{0,K}^2.
\end{eqnarray}
Note that this term actually vanishes for uniformly regular grids (take $\tau_K=h^2$ and integrate by parts).

\medskip

Similarly, for the third term on the r.h.s. of \eqref{eq:ErrEq1}, 
we have:
\begin{eqnarray}\label{eq:cotaIII}
III&=&\sum_{K\in{\cal T}_h}\tau_K (\boldsymbol{\phi}_{r}^{n+1}-\boldsymbol{\phi}_{r}^{n},\nabla s_r^{n+1})_K\nonumber \\
&\leq&  \sum_{K\in{\cal T}_h} \frac{\tau_K}{2\varepsilon\Delta t}\nor{\boldsymbol{\phi}_{r}^{n+1}-\boldsymbol{\phi}_{r}^{n}}{0,K}^2 + \Delta t \sum_{K\in{\cal T}_h}\tau_K \frac{\varepsilon}{2}\nor{\nabla s_r^{n+1}}{0,K}^2\nonumber \\
&\leq& \Delta t\sum_{K\in{\cal T}_h} \frac{C}{\varepsilon}\nor{\boldsymbol{\phi}_{r}^{n+1}-\boldsymbol{\phi}_{r}^{n}}{0,K}^2 + \Delta t \sum_{K\in{\cal T}_h}\tau_K \frac{\varepsilon}{2}\nor{\nabla s_r^{n+1}}{0,K}^2.
\end{eqnarray}

The fourth term on the r.h.s. of \eqref{eq:ErrEq1} can be bounded as follows:
\begin{eqnarray}\label{eq:cotaIV}
IV&=&\Delta t\sum_{K\in{\cal T}_h}\tau_K ( \uv^{n+1}\cdot\nabla\uv^{n+1}-\tilde{\uv}_r^{n+1}\cdot\nabla\uv_r^{n+1},\nabla s_r^{n+1})_K\nonumber \\
&=& \Delta t\sum_{K\in{\cal T}_h}\tau_K (\boldsymbol{\eta}^{n+1}\cdot\nabla\uv^{n+1}-\boldsymbol{\phi}_{r}^{n+1}\cdot\nabla\uv^{n+1},\nabla s_r^{n+1})_K\nonumber \\
&&+\Delta t\sum_{K\in{\cal T}_h}\tau_K (\tilde{\uv}_r^{n+1}\cdot\nabla\boldsymbol{\eta}^{n+1}-\tilde{\uv}_r^{n+1}\cdot\nabla\boldsymbol{\phi}_{r},\nabla s_r^{n+1})_K\nonumber \\
&&+{\Delta t\sum_{K\in{\cal T}_h}\tau_K\left((\uv_r^{n+1}-\tilde{\uv}_r^{n+1})\cdot\nabla\uv^{n+1},\nabla s_r^{n+1}\right)_K}\nonumber \\
&\leq& \Delta t\sum_{K\in{\cal T}_h}\tau_K\nor{\boldsymbol{\eta}^{n+1}}{0,K}\nor{\nabla\uv^{n+1}}{0,\infty,K}\nor{\nabla s_r^{n+1}}{0,K}\nonumber \\
&&+\Delta t\sum_{K\in{\cal T}_h}\tau_K\nor{\boldsymbol{\phi}_{r}^{n+1}}{0,K}\nor{\nabla\uv^{n+1}}{0,\infty,K}\nor{\nabla s_r^{n+1}}{0,K}\nonumber \\
&&+\Delta t\sum_{K\in{\cal T}_h}\tau_K\nor{\tilde{\uv}_r^{n+1}}{0,\infty,K}\nor{\nabla\boldsymbol{\eta}^{n+1}}{0,K}\nor{\nabla s_r^{n+1}}{0,K}\nonumber \\
&&+\Delta t\sum_{K\in{\cal T}_h}\tau_K\nor{\tilde{\uv}_r^{n+1}}{0,\infty,K}\nor{\nabla\boldsymbol{\phi}_{r}^{n+1}}{0,K}\nor{\nabla s_r^{n+1}}{0,K}\nonumber \\
&&+{\Delta t\sum_{K\in{\cal T}_h}\tau_K\nor{\uv_r^{n+1}-\tilde{\uv}_r^{n+1}}{0,K}\nor{\nabla\uv^{n+1}}{0,\infty,K}\nor{\nabla s_r^{n+1}}{0,K}}\nonumber \\
&\leq&\Delta t\sum_{K\in{\cal T}_h}\frac{C}{\varepsilon}\tau_K(\nor{\boldsymbol{\eta}^{n+1}}{0,K}^2+\nor{\boldsymbol{\phi}_{r}^{n+1}}{0,K}^2)+5\Delta t\sum_{K\in{\cal T}_h}\frac{\varepsilon}{2}\tau_K\nor{\nabla s_r^{n+1}}{0,K}^2\nonumber \\
&&+\Delta t\sum_{K\in{\cal T}_h}\frac{C}{\varepsilon}\frac{\tau_K}{h_K^2}\left(\nor{\nabla\boldsymbol{\eta}^{n+1}}{0,K}^2+\nor{\nabla\boldsymbol{\phi}_{r}^{n+1}}{0,K}^2\right)\nonumber \\
&&+{\Delta t\sum_{K\in{\cal T}_h}\frac{C}{\varepsilon}\tau_K(\nor{\uv_r^{n+1}}{0,K}^2+\nor{\tilde{\uv}_r^{n+1}}{0,K}^2)}\nonumber \\
&\leq&\Delta t\sum_{K\in{\cal T}_h}\frac{C}{\varepsilon}\tau_K(\nor{\boldsymbol{\eta}^{n+1}}{0,K}^2+\nor{\boldsymbol{\phi}_{r}^{n+1}}{0,K}^2)+5\Delta t\sum_{K\in{\cal T}_h}\frac{\varepsilon}{2}\tau_K\nor{\nabla s_r^{n+1}}{0,K}^2\nonumber \\
&&+\Delta t\sum_{K\in{\cal T}_h}\frac{C}{\varepsilon}\frac{\tau_K}{h_K^2}\left(\nor{\nabla\boldsymbol{\eta}^{n+1}}{0,K}^2+\nor{\nabla\boldsymbol{\phi}_{r}^{n+1}}{0,K}^2\right)+{\Delta t\frac{C}{\varepsilon}h^2},
\end{eqnarray} 
where we used H\"older's, Young's and triangle inequalities, the velocity stability estimate $\eqref{estvel1}$, and \eqref{hip:utrunc} to bound $\tilde{\uv}_r^{n+1}$ when $d=3$. {Note that the last penalty term in \eqref{eq:cotaIV} vanishes when $d=2$, since $\tilde{\uv}_r^{n+1}=\uv_r^{n+1}$}.

\medskip

Using local inverse estimates \cite{Bernardi04}, we can bound the fifth term on the r.h.s. of \eqref{eq:ErrEq1} as:
\begin{eqnarray}\label{eq:cotaV}
V&=&-\nu\Delta t\sum_{K\in{\cal T}_h}\tau_K \left( \Delta (\boldsymbol{\eta}^{n+1}-\boldsymbol{\phi}_r^{n+1}),\nabla s_r^{n+1}\right)_K\nonumber \\
&\leq& \nu^2\Delta t\sum_{K\in{\cal T}_h}\frac{\tau_K}{2\varepsilon}\left(\nor{\Delta \boldsymbol{\eta}^{n+1}}{0,K}^2+\nor{\Delta \boldsymbol{\phi}_r^{n+1}}{0,K}^2\right)+2 \Delta t\sum_{K\in{\cal T}_h}\frac{\varepsilon}{2}\tau_K\nor{\nabla s_r^{n+1}}{0,K}^2\nonumber \\
&\leq& \frac{C}{\varepsilon}\nu^2\,\Delta t(\nor{\Delta\boldsymbol{\eta}^{n+1}}{\tau}^2
+\nor{\nabla\boldsymbol{\phi}_r^{n+1}}{0}^2)
+2 \Delta t\sum_{K\in{\cal T}_h}\frac{\varepsilon}{2}\tau_K\nor{\nabla s_r^{n+1}}{0,K}^2,
\end{eqnarray}

\medskip

Finally, by Cauchy-Schwarz and Young's inequalities, we bound the sixth term on the r.h.s. of \eqref{eq:ErrEq1}:
\begin{eqnarray}\label{eq:cotaVI}
VI=\Delta t\sum_{K\in{\cal T}_h}\tau_K ( \nabla \rho^{n+1},\nabla s_r^{n+1})_K&\leq& \Delta t \sum_{K\in{\cal T}_h} \frac{\tau_K}{2\varepsilon}\nor{\nabla \rho^{n+1}}{0,K}^2 \nonumber \\
&&+ \Delta t \sum_{K\in{\cal T}_h}\tau_K \frac{\varepsilon}{2}\nor{\nabla s_r^{n+1}}{0,K}^2.
\end{eqnarray}

Summarizing \eqref{eq:cotaI}-\eqref{eq:cotaVI} and taking $\varepsilon=1/11$ from \eqref{eq:ErrEq1}, we get:
\begin{eqnarray}\label{eq:ErrEq2}
\Delta t\sum_{K\in{\cal T}_h}\tau_K\nor{\nabla s_r^{n+1}}{0,K}^2&\leq& C\Delta t\left(h^2\nor{\boldsymbol{c}^n}{0}^2+\nor{\boldsymbol{\eta}^{n+1}-\boldsymbol{\eta}^{n}}{0}^2+\nor{\boldsymbol{\phi}_r^{n+1}-\boldsymbol{\phi}_r^{n}}{0}^2\right)\nonumber \\
&&+C\Delta t h^2(\nor{\boldsymbol{\eta}^{n+1}}{0}^2+\nor{\boldsymbol{\phi}_r^{n+1}}{0}^2+{1})
\nonumber \\
&&+C\Delta t (\nor{\nabla\boldsymbol{\eta}^{n+1}}{0}^2+\nor{\nabla\boldsymbol{\phi}_r^{n+1}}{0}^2)\nonumber \\
&&+C\nu^2\Delta t\nor{\Delta\boldsymbol{\eta}^{n+1}}{\tau}^2+C\Delta t\nor{\nabla\rho^{n+1}}{\tau}^2.
\end{eqnarray}

Summing \eqref{eq:ErrEq2} from $n=0$ to $N-1$, we have:
\begin{eqnarray}\label{eq:ErrEq3}
\Delta t\sum_{n=0}^{N-1}\nor{\nabla s_r^{n+1}}{\tau}^2&\leq&C\Delta t\sum_{n=0}^{N-1}\left(h^2 \nor{\boldsymbol{c}^n}{0}^2+\nor{\boldsymbol{\eta}^{n+1}-\boldsymbol{\eta}^{n}}{0}^2
+\nor{\boldsymbol{\phi}_r^{n+1}-\boldsymbol{\phi}_r^{n}}{0}^2\right)\nonumber \\ 
&&+C\Delta t \sum_{n=0}^{N-1}h^2(\nor{\boldsymbol{\eta}^{n+1}}{0}^2+\nor{\boldsymbol{\phi}_r^{n+1}}{0}^2)+{C h^2}\nonumber \\
&&+C\Delta t \sum_{n=0}^{N-1}(\nor{\nabla\boldsymbol{\eta}^{n+1}}{0}^2+\nor{\nabla\boldsymbol{\phi}_r^{n+1}}{0}^2)\nonumber \\
&&+C\Delta t\sum_{n=0}^{N-1}\left(\nu^2\nor{\Delta\boldsymbol{\eta}^{n+1}}{\tau}^2+\nor{\nabla\rho^{n+1}}{\tau}^2\right).
\end{eqnarray}
Next, we estimate each term on the r.h.s. of \eqref{eq:ErrEq3}.
%

\medskip

By Taylor's theorem, the first term on the r.h.s. of \eqref{eq:ErrEq3} can be estimated as follows:
\BEQ\label{eq:cn}
\Delta t\sum_{n=0}^{N-1}h^2\nor{\boldsymbol{c}^{n}}{{\bf L}^{2}}^{2}\leq C\Delta t^{2}h^2\nor{\partial_{tt}\uv}{L^{2}({\bf L}^{2})}^{2}\leq C\Delta t^{2}h^2.
\EEQ

\medskip

To estimate the second term on the r.h.s. of \eqref{eq:ErrEq3}, we use triangle inequality and the $L^2$ error estimate for $P_r^v \uv$ (see \cite{IliescuWang14}, Lemma 3.3):
\BEQ\label{eq:eta2}
\Delta t\sum_{n=0}^{N-1}\nor{\boldsymbol{\eta}^{n+1}-\boldsymbol{\eta}^{n}}{0}^2\leq C\left(h^{2l}+\Delta t^2+\sum_{i=r+1}^{M_v}\lambda_i\right).
\EEQ

By using triangle and Cauchy-Schwarz inequalities, the bound \eqref{eq:cota_grad_div}, and estimate (73) in \cite{NovoRubinoSINUM21}, the third term on the r.h.s. of \eqref{eq:ErrEq3} can be estimated as follows:
\begin{eqnarray}\label{eq:phi2}
\Delta t\sum_{n=0}^{N-1}\nor{\boldsymbol{\phi}_r^{n+1}-\boldsymbol{\phi}_r^{n}}{0}^2 &\leq& 2 \Delta t\sum_{n=0}^{N-1}\left(\nor{(\uv_r^{n+1}-P_r^{v}\uv_h^{n+1})-(\uv_r^{n}-P_r^{v}\uv_h^{n})}{0}^2\right)\nonumber \\
&&+2 \Delta t\sum_{n=0}^{N-1}\left(\nor{P_r^v\left[(\uv_h^{n+1}-\uv^{n+1})-(\uv_h^{n}-\uv^{n})\right]}{0}^2\right)\nonumber \\
&\leq& C\sum_{i=r+1}^{M_v}\lambda_i\nor{\nabla \boldsymbol{\varphi}_{i}}{0}^2\nonumber \\
&&+C \Delta t\sum_{n=0}^{N-1}\left(\nor{(\uv_h^{n+1}-\uv^{n+1})-(\uv_h^{n}-\uv^{n})}{0}^2\right)\nonumber \\
&\leq&C\left(\sum_{i=r+1}^{M_v}\lambda_i\nor{\nabla \boldsymbol{\varphi}_{i}}{0}^2 + h^{2l}+\Delta t^2\right).
\end{eqnarray}

\medskip

Similarly, for the fourth and fifth term on the r.h.s. of \eqref{eq:ErrEq3} we have:
\BEQ\label{eq:eta1}
\Delta t\sum_{n=0}^{N-1}h^2\nor{\boldsymbol{\eta}^{n+1}}{0}^2\leq C\,h^2\left(h^{2l}+\Delta t^2+\sum_{i=r+1}^{M_v}\lambda_i\right),
\EEQ
and
\begin{eqnarray}\label{eq:phi1}
\Delta t\sum_{n=0}^{N-1}h^2\nor{\boldsymbol{\phi}_r^{n+1}}{0}^2 &\leq& 2 \Delta t\sum_{n=0}^{N-1}h^2\left(\nor{\uv_r^{n+1}-P_r^{v}\uv_h^{n+1}}{0}^2+\nor{P_r^v(\uv_h^{n+1}-\uv^{n+1})}{0}^2\right)\nonumber \\
&\leq&C\,h^2\left(\sum_{i=r+1}^{M_v}\lambda_i\nor{\nabla \boldsymbol{\varphi}_{i}}{0}^2 + h^{2l}+\Delta t^2\right).
\end{eqnarray}

\medskip

For the sixth term on the r.h.s. of \eqref{eq:ErrEq3}, we have:
\BEQ\label{eq:eta1grad}
\Delta t\sum_{n=0}^{N-1}\nor{\nabla\boldsymbol{\eta}^{n+1}}{0}^2
\leq C\left[(h^{2l}+\Delta t^2)(\nu^{-1}+\nor{S_r^v}{2})+\sum_{i=r+1}^{M_v}\lambda_i\nor{\nabla \boldsymbol{\varphi}_{i}}{0}^2\right],
\EEQ
where we used the $H_0^1$ error estimate for $P_r^v \uv$ (see \cite{IliescuWang14}, Lemma 3.3), and $\nor{S_r^v}{2}$ denotes the $2$-norm of the stiffness velocity matrix with entries $[S_r^v]_{ij}=(\nabla\boldsymbol{\varphi}_i,\nabla\boldsymbol{\varphi}_j)$, $i,j=1,\ldots,r$ that comes from the use of POD inverse estimates (see \cite{KunischVolkwein01}, Lemma 2).

\medskip

For the seventh term on the r.h.s. of \eqref{eq:ErrEq3}, we have:
\begin{eqnarray}\label{eq:phi1grad}
\Delta t\sum_{n=0}^{N-1}\nor{\nabla\boldsymbol{\phi}_r^{n+1}}{0}^2 &\leq& 2\Delta t \sum_{n=0}^{N-1}\left(\nor{\nabla(\uv_r^{n+1}-P_r^{v}\uv_h^{n+1})}{0}^2+\nor{\nabla\left[P_r^v(\uv_h^{n+1}-\uv^{n+1})\right]}{0}^2\right)\nonumber \\
&\leq& C\nor{S_r^v}{2}\left(\sum_{i=r+1}^{M_v}\lambda_i\nor{\nabla \boldsymbol{\varphi}_{i}}{0}^2 + h^{2l}+\Delta t^2\right),
\end{eqnarray}
where we used the triangle inequality, POD inverse estimates (see \cite{KunischVolkwein01}, Lemma 2), estimate (73) in \cite{NovoRubinoSINUM21}, and bound \eqref{eq:cota_grad_div}.

\medskip

The eight term on the r.h.s. of \eqref{eq:ErrEq3} can be bounded as follows:
\begin{eqnarray}\label{eq:eta1lap}
\Delta t\sum_{n=0}^{N-1}\nu^2\nor{\Delta\boldsymbol{\eta}^{n+1}}{\tau}^2
&\leq& C\nu^2\Delta t\sum_{n=0}^{N-1}\left[\nor{\Delta(\uv^{n+1}-\uv_h^{n+1})}{\tau}^2+\nor{\Delta(\uv_h^{n+1}-P_r^v \uv_h^{n+1})}{\tau}^2\right]\nonumber \\
&&+C\nu^2\Delta t\sum_{n=0}^{N-1}\nor{\Delta[P_r^v(\uv_h^{n+1}-\uv^{n+1})]}{\tau}^2\nonumber \\
&\leq& C\nu^2\Delta t\sum_{n=0}^{N-1}\left(h^{2l}+\nor{\nabla(\uv^{n+1}-\uv_h^{n+1})}{0}^2+\nor{\nabla(\uv_h^{n+1}-P_r^v \uv_h^{n+1})}{0}^2\right)\nonumber \\
&&+C\nu^2\nor{S_r^v}{2}\,\Delta t\sum_{n=0}^{N-1}\nor{\uv^{n+1}-\uv_h^{n+1}}{0}^2\nonumber \\
&\leq& C\left[\nu^2\,h^{2l}+\nu(1+\nu\nor{S_r^v}{2})(h^{2l}+\Delta t^2)+\nu^2\sum_{i=r+1}^{M_v}\lambda_i\nor{\nabla\boldsymbol{\varphi}_i}{0}^2\right],
\end{eqnarray}
where we used the triangle inequality, optimal and local inverse estimates \cite{Bernardi04}, POD inverse estimates (see \cite{KunischVolkwein01}, Lemma 2), bound \eqref{eq:cota_grad_div}, and the POD projection error in $H^1$-seminorm for $\uv_h$ (see \cite{IliescuJohn15}, Lemma 3.2).

\medskip

Finally, for the last term on the r.h.s. of \eqref{eq:ErrEq3}, similarly to \eqref{eq:eta1grad} we obtain:
\begin{eqnarray}\label{eq:rho1grad}
\Delta t\sum_{n=0}^{N-1}\nor{\nabla\rho^{n+1}}{\tau}^2
&\leq& C\Delta t\sum_{n=0}^{N-1}\left[\nor{\nabla(p^{n+1}-p_h^{n+1})}{\tau}^2+\nor{\nabla(p_h^{n+1}-P_r^p p_h^{n+1})}{\tau}^2\right]\nonumber \\
&&+C\Delta t\sum_{n=0}^{N-1}\nor{\nabla[P_r^p(p_h^{n+1}-p^{n+1})]}{\tau}^2\nonumber \\
&\leq&C\Delta t\sum_{n=0}^{N-1}\left(h^{2l}+\nor{p^{n+1}-p_h^{n+1}}{0}^2+h^2\nor{\nabla(p_h^{n+1}-P_r^p p_h^{n+1})}{0}^2\right)\nonumber \\
&\leq& C\left[h^{2l}+h^{-1}(h^{2l}+\Delta t^2)+h^2\sum_{i=r+1}^{M_p}\gamma_i\nor{\nabla\psi_i}{0}^2\right],
\end{eqnarray}
where we used the triangle inequality, optimal and local inverse estimates \cite{Bernardi04}, bound \eqref{eq:cota_grad_div_pre}, and the POD projection error in $H^1$-seminorm for $p_h$ (see \cite{IliescuJohn15}, Lemma 3.2).
%

\medskip

Collecting \eqref{eq:cn}-\eqref{eq:rho1grad}, estimate \eqref{eq:ErrEq3} becomes:
\begin{eqnarray}\label{eq:ErrEq4}
\Delta t\sum_{n=0}^{N-1}\nor{\nabla s_r^{n+1}}{\tau}^2&\leq&C\left[{h^2}+(h^{2l}+\Delta t^2)(1+\nu^{-1}+\nor{S_r^v}{2}+h^{-1})\right]\nonumber \\
&& +C\sum_{i=r+1}^{M_v}\lambda_i\left[1+\left(1+\nor{S_r^v}{2}\right)\nor{\nabla\boldsymbol{\varphi}_i}{0}^2\right]\nonumber \\
&&+C\,h^2\sum_{i=r+1}^{M_p}\gamma_i\nor{\nabla\psi_i}{0}^2.
\end{eqnarray}

\medskip

Finally, using the triangle inequality and estimate \eqref{eq:rho1grad}, we conclude:
\begin{eqnarray}\label{eq:ErrEq5}
\Delta t\sum_{n=0}^{N-1}\nor{\nabla(p_r^{n+1}-p^{n+1})}{\tau}^2&\leq& 2\Delta t\sum_{n=0}^{N-1}\left(\nor{\nabla\rho^{n+1}}{\tau}^2+\nor{\nabla s_r^{n+1}}{\tau}^2\right) \nonumber \\
&\leq&C\left[{h^2}+(h^{2l}+\Delta t^2)(1+\nu^{-1}+\nor{S_r^v}{2}+h^{-1})\right]\nonumber \\
&& +C\sum_{i=r+1}^{M_v}\lambda_i\left[1+\left(1+\nor{S_r^v}{2}\right)\nor{\nabla\boldsymbol{\varphi}_i}{0}^2\right]\nonumber \\
&&+C\,h^2\sum_{i=r+1}^{M_p}\gamma_i\nor{\nabla\psi_i}{0}^2,
\end{eqnarray}
{where we recall that the penalty term $h^2$ in \eqref{eq:ErrEq5} vanishes when $d=2$}. This concludes the proof.
\qed

\begin{remark}\label{rm:AltPresEst}
An alternative pressure error estimate can be obtained by using a different bound for the term $IV$ in \eqref{eq:cotaIV} when $d=3$, taking advantage of the fact that $\tilde{\uv}_r^{n+1}= {\uv}_r^{n+1}$ when  $\nor{{\uv}_r^{n+1}}{0,\infty,K} \le C h_K^{-1}$. Indeed, we can replace in \eqref{eq:cotaIV} the last term by:
\begin{eqnarray}\label{eq:cotaIValt}
{\Delta t\sum_{K\in{\cal T}_h^*}\frac{C}{\varepsilon}\tau_K(\nor{\uv_r^{n+1}}{0,K}^2+\nor{\tilde{\uv}_r^{n+1}}{0,K}^2)}&\leq&{\Delta t\sum_{K\in{\cal T}_h^*}\frac{C}{\varepsilon}\tau_K |K|(\nor{\uv_r^{n+1}}{0,\infty,K}^2+\nor{\tilde{\uv}_r^{n+1}}{0,\infty,K}^2)}\nonumber \\
&\leq&{\Delta t\sum_{K\in{\cal T}_h^*}\frac{C}{\varepsilon}h_K^2 |K|\left(h_K^{-3}\nor{\uv_r^{n+1}}{0,K}^2+h_K^{-2}\right)}\nonumber \\
&\leq& {\Delta t\,\frac{C}{\varepsilon}\,h^{-1}|\Om^*|},
\end{eqnarray}
where ${\cal T}_h^*=\{K\in{\cal T}_h : \nor{\uv_r^{n+1}}{0,\infty,K}\geq C h_K^{-1}\}$, $|\Om^*|$ is the measure of ${\cal T}_h^*$ (i.e., $|\Om^*|=\disp\sum_{K\in{\cal T}_h^*}|K|$), and we used local inverse estimates \cite{Bernardi04}, the velocity stability estimate $\eqref{estvel1}$, and \eqref{hip:utrunc} to bound $\tilde{\uv}_r^{n+1}$, assuming also uniformly regular grids. 

\medskip

Now, we have:
\begin{eqnarray}\label{eq:cotaOmStar}
|\Om^*|&=&\sum_{K\in{\cal T}_h^*}\int_K \frac{h_K}{h_K}\,d\xv\leq C\sum_{K\in{\cal T}_h^*}\int_K h_K \nor{\uv_r^{n+1}}{0,\infty,K}d\xv\nonumber \\
&\leq& C\sum_{K\in{\cal T}_h^*}\int_K h_K^{1/2}\nor{\uv_r^{n+1}}{0,6,K}d\xv\nonumber 
\leq C\sum_{K\in{\cal T}_h^*}\nor{\nabla\uv_r^{n+1}}{0,K}\,|K|\,h_K^{1/2}\nonumber \\
&\leq& C\left(\sum_{K\in{\cal T}_h^*} \nor{\nabla\uv_r^{n+1}}{0,K}^2\right)^{1/2}\left(\sum_{K\in{\cal T}_h^*}|K|^2 h_K\right)^{1/2}\nonumber \\
&\leq& C\nor{\nabla\uv_r^{n+1}}{0,\Omega}\left(\sum_{K\in{\cal T}_h^*}|K| h_K^4\right)^{1/2}\leq C\nor{\nabla\uv_r^{n+1}}{0,\Omega}h^2|\Om^*|^{1/2}.
\end{eqnarray}
where we used the definition of ${\cal T}_h^*$, local inverse estimates \cite{Bernardi04}, H\"older's inequality, the Sobolev embedding $H^1\hookrightarrow L^6$, and Cauchy-Schwarz inequality. Thus, from \eqref{eq:cotaOmStar} we get:
\BEQ\label{eq:CotaOmStarDef}
|\Om^*|\leq C\nor{\nabla\uv_r^{n+1}}{0,\Omega}^2 \,h^4, 
\EEQ
and substituting \eqref{eq:CotaOmStarDef} in \eqref{eq:cotaIValt} we obtain:
\BEQ\label{eq:cotaIValtDef}
{\Delta t\sum_{K\in{\cal T}_h^*}\frac{C}{\varepsilon}\tau_K(\nor{\uv_r^{n+1}}{0,K}^2+\nor{\tilde{\uv}_r^{n+1}}{0,K}^2)}\leq {\Delta t\,\frac{C}{\varepsilon}\,h^3 \nor{\nabla\uv_r^{n+1}}{0,\Omega}^2}.
\EEQ

\medskip

Thus, by applying the velocity stability estimate $\eqref{estvel2}$, in this case we can replace the penalty term $h^2$ in \eqref{eq:Estar} by $h^3\nu^{-1}$.
\end{remark}

\section{Numerical tests} \label{sec:NumStud}
The data used for the construction of the POD subspaces are obtained by using Fenics \cite{FenicsBook} with Taylor-Hood element $\mbb{P}_2 / \mbb{P}_1 $. 
To quantify the error estimates derived in the previous section, we consider two classical benchmark numerical examples, the flow past a cylinder and the lid-driven cavity flow. Since the errors depend strongly on the velocity POD truncation order, the error analyses are presented in terms of the velocity POD contribution ratio given as
\begin{equation}\label{ratio_POD}
\mathcal{R}_{u} = 1-\somme{i}{1}{r}\lambda_i / \somme{i}{1}{\mcal{N}}\lambda_i\review{,}
\end{equation}
\ModifMourad{where $r$ and $\mcal{N}$ are respectively the truncation order and the total number of snapshots. This} contribution ratio can be seen as a refinement parameter of the error. It can be used along with the predicted error-slope to determine a priori the accuracy of the SM-ROM.
 
\subsection{Flow past a cylinder Re = 200}
\begin{figure}[hbtp!]
\centering 
\includegraphics[width=0.85\linewidth]{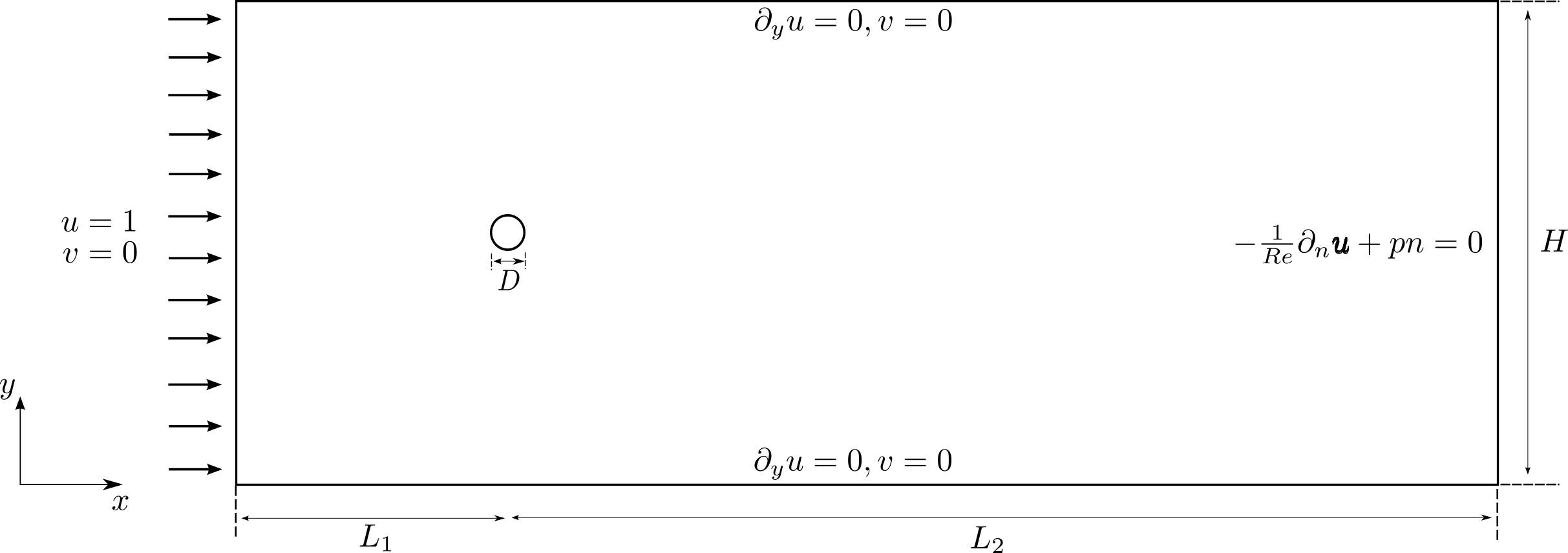}
\caption{Two-dimensional domain and boundary conditions for the problem of flow past a circular cylinder.}
\label{Fig:Flow_past_cylinder_bd_conditions}
\end{figure}

Consider the two dimensional flow past a cylinder of diameter $D$, placed in a rectangular domain of length $H = 30D$ and width = $45 D$. The center of the cylinder is situated at $L_1 = 10D$ from the inlet and $H/2$ from horizontal walls. The flow is driven by an inlet velocity $U$ and is allowed to flow through the outlet. Free slip boundary conditions are applied at the horizontal edges while no slip boundary condition are considered on the cylinder's wall. Illustration of boundary condition is given in Figure \ref{Fig:Flow_past_cylinder_bd_conditions}.

\begin{figure}[hbtp!]
\centering 
\includegraphics[width=0.5\linewidth]{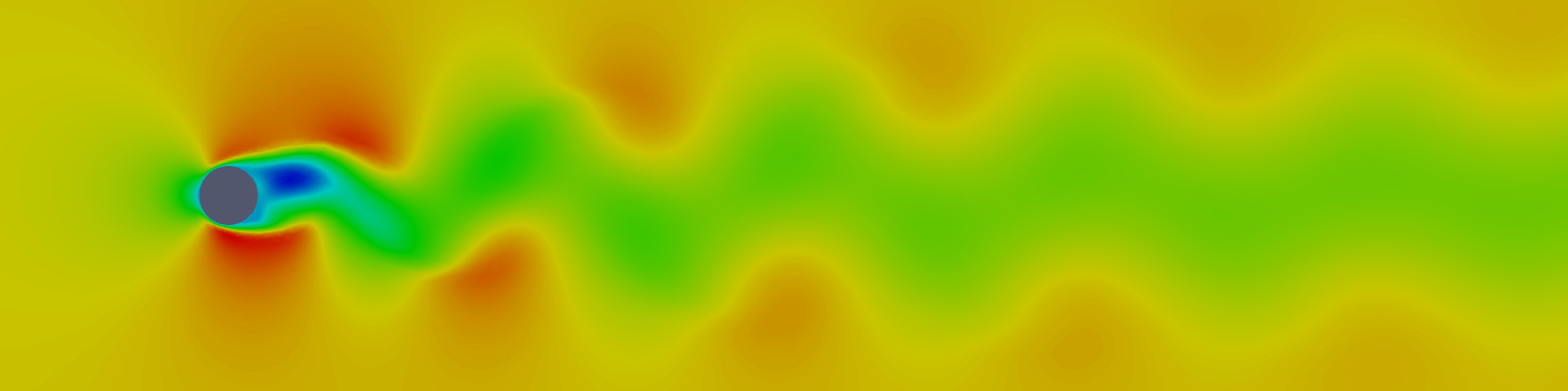}
\\\vspace{0.3cm}
\includegraphics[width=0.5\linewidth]{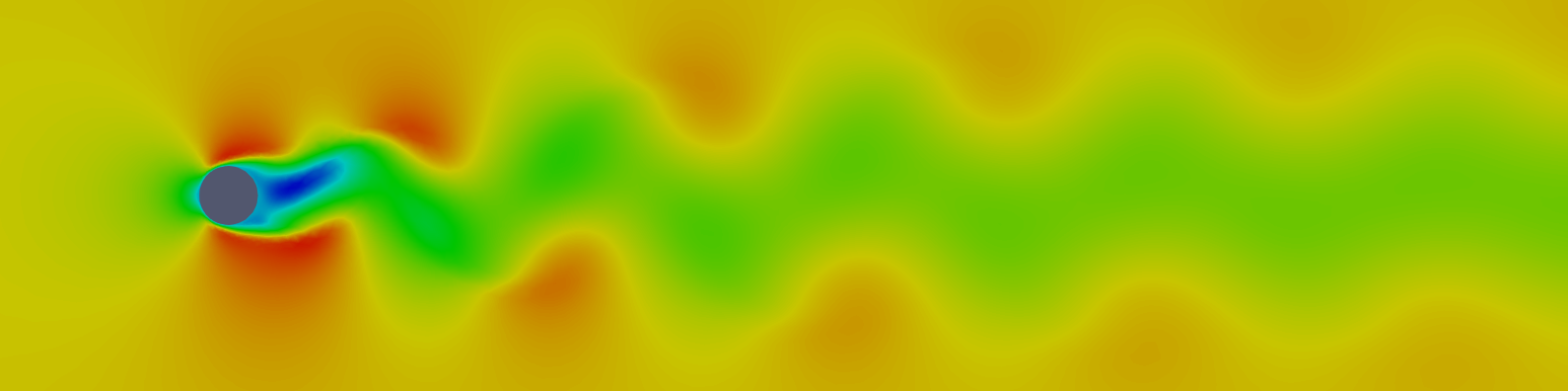}
\\\vspace{0.3cm}
\includegraphics[width=0.5\linewidth]{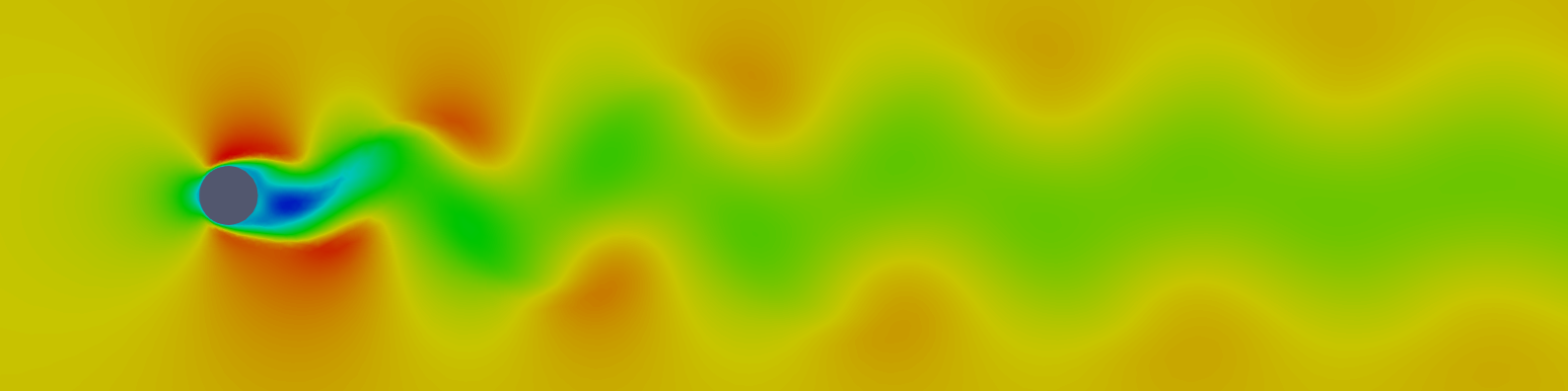}
\caption{Illustration of Von Karman vortices arising in the flow past a cylinder for $Re = 200$ at \ModifMourad{three different time instants.}}
\label{Fig:flow past cylinder}
\end{figure}

\begin{figure}[hbtp!]
\centering 
\includegraphics[width=\linewidth]{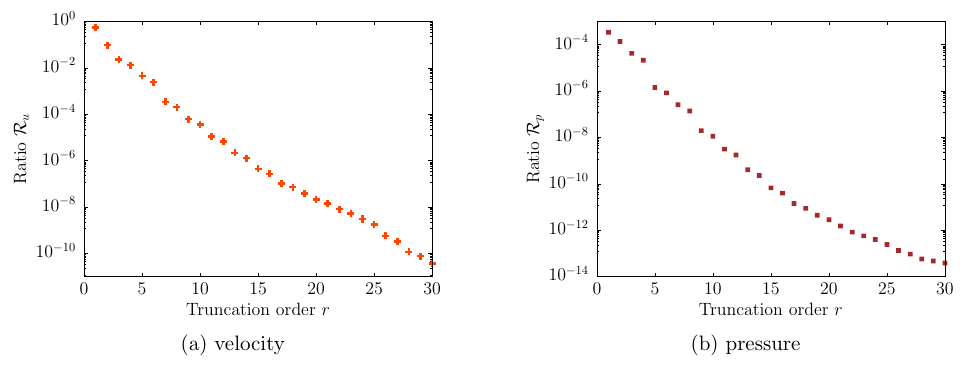}
\caption{Evolution of the velocity and pressure POD contribution ratio with respect to modes numbers for the flow past a cylinder at $Re=200$.}
\label{RC_RE150}
\end{figure}

For the numerical simulation, the semi-implicit Euler method is used with time step $\Delta t= 10^{-2} $, while a \ModifMourad{non-uniform triangular mesh made of $21174$ cells} is used for the spatial discretization. At $Re = 200$, the flow creates alternating low-high pressure vortices downstream the cylinder,  triggering the generation of the Von Karman vortex pattern in the wake of the flow past the cylinder. These structures are illustrated in Figure \ref{Fig:flow past cylinder}. 
\\
In order to construct the velocity and pressure POD bases, \ModifMourad{$400$ uniformly distributed snapshots covering $8$ periods of the periodic regime of the flow are considered}. The evolution of the contribution ratio of POD modes is given as a function of the modes numbers in Figure \ref{RC_RE150}. A first inspection of the contribution ratio reveals that the first three modes carry the biggest share of the flow energy. After that, the contribution evolves by couple of modes. This behavior is essentially due the double multiplicity of POD eigenvalues revealing the presence of eigenspaces of dimension two. In that case, it's worth noting that the SM-ROM error is most likely to oscillate if the POD basis is truncated at a mode forming the first vector of a two dimensional POD eigenspace. To prevent that from happening, it is preferred to construct the SM-ROM gradually with respect to eigenspaces rather than to eigenvectors. 
%

\begin{figure}[hbtp!]
\centering 
\includegraphics[width=\linewidth]{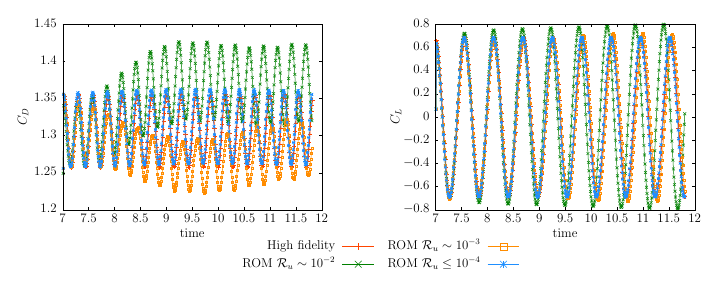}
\caption{Hydrodynamics coefficients}
\label{Hydrodynamics coefficients}
\end{figure}

\begin{figure}[hbtp!]
\centering 
\includegraphics[width=\linewidth]{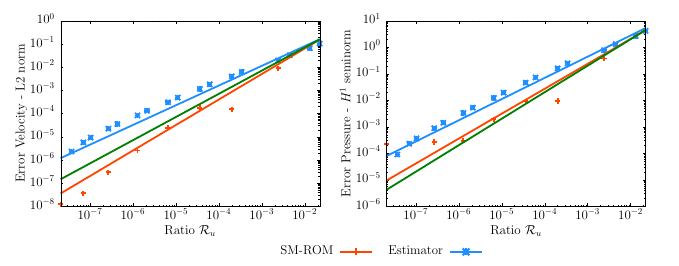}
\caption{\ModifMourad{Log-Log scale representation of the estimator versus SM-ROM errors for the flow past a cylinder at $Re=200$}. The dots refer to the calculated values while the solid lines are their power regression. The green solid line corresponds to the curve of logarithmic slope equal to $1$. \ModifMourad{Only the errors corresponding to $\mcal{R}_u$ larger than $10^{-6}$ have been taken into account to calculate the regression lines.}}
\label{CD}
\label{ErrEstim_RE200}
\end{figure}

The POD contribution ratio in Figure \ref{RC_RE150} tends rapidly to the numerical zero meaning by that the few first modes are sufficiently enough to catch the quasi-totality of the flow energy. 
The prediction quality of the SM-ROM improves with decreasing values of $\mathcal{R}_u$. This can be seen from the hydrodynamics coefficients $C_D$ and $C_L$ given by 
\begin{equation*}
C_D = \myfrac{F_D}{\frac{1}{2}\rho U^2 D} \ \ , \ \  C_L = \myfrac{F_L}{\frac{1}{2}\rho U^2 D}\review{,}
\end{equation*}
\ModifMourad{where $\rho$ and $U$ are respectively the density of the fluid and the characteristic velocity}, and $F_D$ and $F_L$ are the drag and lift forces exerted by the fluid on the cylinder's boundary
\begin{equation*}
\left(\begin{matrix}
F_D \\
F_L
\end{matrix}\right)  =  \int_{{_{\txt{cylinder}}}} \left( \mu (\nabla \vect{u} + \nabla^T \vect{u}) - p I \right)\cdot \vect{n} \,d\sigma\review{,}
\end{equation*} 
\ModifMourad{Here $\mu$ refers to the dynamic viscosity}. Figure \ref{Hydrodynamics coefficients} gives a confrontation of the hydrodynamic coefficients for variable magnitudes of the ratio $\mathcal{R}_u$. It shows that the predicted coefficients by SM-ROM converge towards the reference coefficients as $\mathcal{R}_u$ drops down to zero. Precisely, when the ratio attains an order of $10^{-4}$ or less, the solutions are well recovered and the SM-ROM coefficients are in good agreement with the underlying high fidelity model. 
\\
In the following, we investigate the errors evolution with respect to the POD contribution ratio.
In theory, the error by the SM-ROM and the POD contribution ratio must decrease simultaneously. However, this is not always the case. It can be noticed from Figure \ref{ErrEstim_RE200} that while the SM-ROM velocity errors decrease as expected, the SM-ROM pressure error experiences a stagnation for ratios smaller than $10^{-6}$, meaning that the SM-ROM has attained its maximum precision at this stage.  {This can be explained as the high fidelity weak formulation \eqref{eq:discTempFOM} is not exactly verified for gradients of pressure POD modes, because these do not lie in the discrete velocity space $\bm{X}_h$. This effect does not appear in the error estimates \eqref{eq:cota_finalSUPp}, as these apply to the error between POD-ROM and exact pressures. This limit accuracy for the reduced pressure computation should decrease as the discretization parameters $h$ and $\Delta t$ decrease. }
%
\\
Reduced order model errors with respect to the rate of convergence are reported in Figure \ref{ErrEstim_RE200}.
For the sake of clarity, and since we are interested only in the error slopes, the constant in the error estimation is adjusted in such a way that the estimated and SM-ROM errors start from the same point. Moreover, a power law regression with solid lines is overlayed on the computed errors and a curve of logarithmic slope equal to one is plotted alongside for comparison. Given the likely stagnation behavior of the errors at some values of the contribution ratio, the regression is only applied to the first points where variations are noticed. 

It can be seen that the logarithmic slops for the SM-ROM errors are close to one meaning that the rate of convergence of the errors with respect to $\mathcal{R}_u$ is close to $1$. In other words, a jump of one order of magnitude in the ratio $\mathcal{R}_u$ results in a jump of the same order in the error. This observation is in good agreement with Figure \ref{Hydrodynamics coefficients} which shows that 
\ModifMourad{
between POD ratios $10^{-2}$ and $10^{-4}$, the SM-ROM jumps two orders of magnitudes and attains $4\times 10^{-4}$ for velocity and $2\times 10^{-2}$ for pressure.
}
These errors represent a sufficient accuracy to reproduce the dynamics of the flow and thus meet the high fidelity model hydrodynamic coefficients. 

Now by inspecting the estimator, one can see that even though the slopes are slightly smaller, it provides a good indicator of the convergence rate of the SM-ROM. If we look at the regression of the estimated errors at ratio $10^{-4}$, it indicates a value of \ModifMourad{$\frac{3}{2}\times 10^{-3}$} for velocity and \ModifMourad{$6\times 10^{-2}$} for the pressure. These are sufficiently good approximations of the actual SM-ROM error that confirm the robustness of the proposed estimator. 


\subsection{Test case 2 : Lid driven cavity Re = 9500}
In this test case we assess the error estimates on the two dimensional lid driven cavity flow problem with $Re=9500$. The problem domain consists of a square cavity $]0,D[\times]0,D[$ where the fluid is driven by a tangential velocity $U$ applied at the top wall. The remaining walls are defined as no-slip conditions. The high fidelity problem is solved in \ModifMourad{a triangular mesh composed of $32928$ cells} by using the semi-implicit Euler method with a fixed time step $\Delta t= 10^{-3} $.

\begin{figure}[hbtp!]
\centering 
\includegraphics[width=0.3\linewidth]{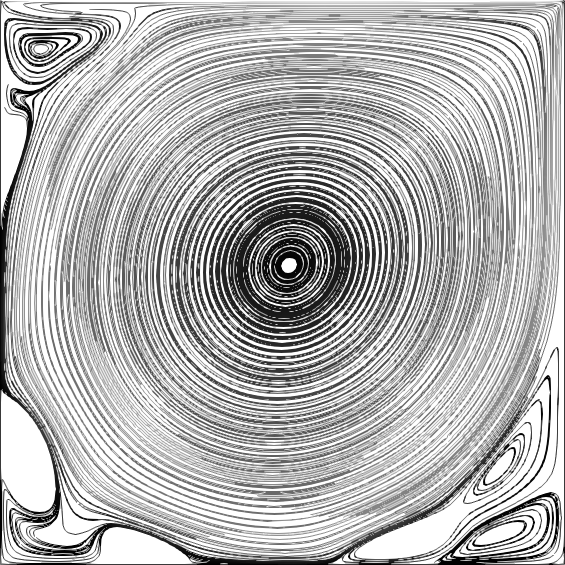}%
\quad\includegraphics[width=0.3\linewidth]{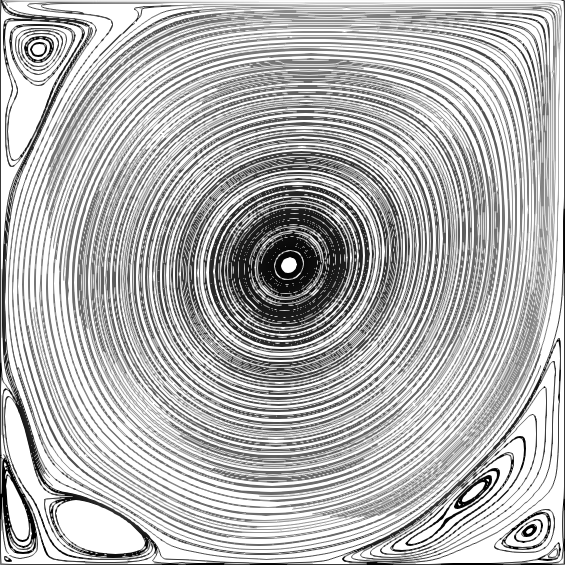}%
\quad\includegraphics[width=0.3\linewidth]{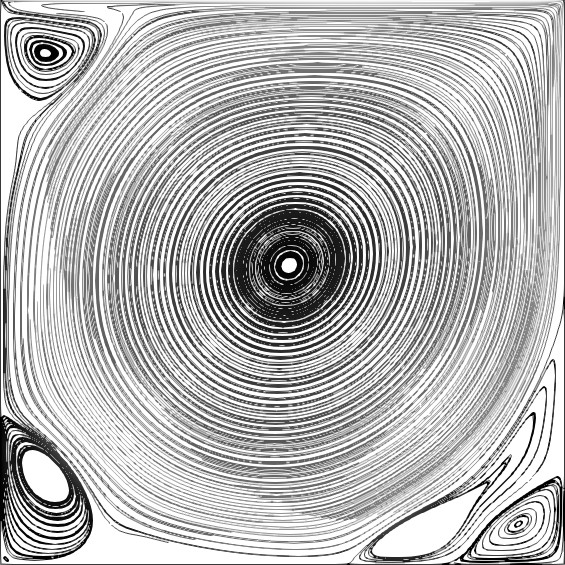}
\caption{\ModifMourad{Streamlines in three different instants of the periodic regime of the flow in a lid driven cavity for $Re=9500$.}}
\label{Fig:flow cavity}
\end{figure}

\begin{figure}[hbtp!]
\centering 
\includegraphics[width=\linewidth]{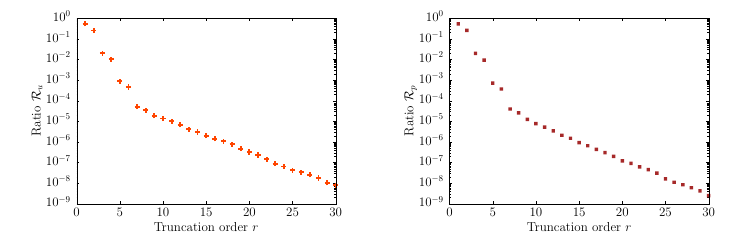}
\caption{Evolution of the velocity and pressure POD contribution ratio with respect to modes numbers for the Lid driven cavity problem at $Re=9500$.}
\label{RC_RE9500}
\end{figure}

\begin{figure}[hbtp!]
\centering 
\includegraphics[width=\linewidth]{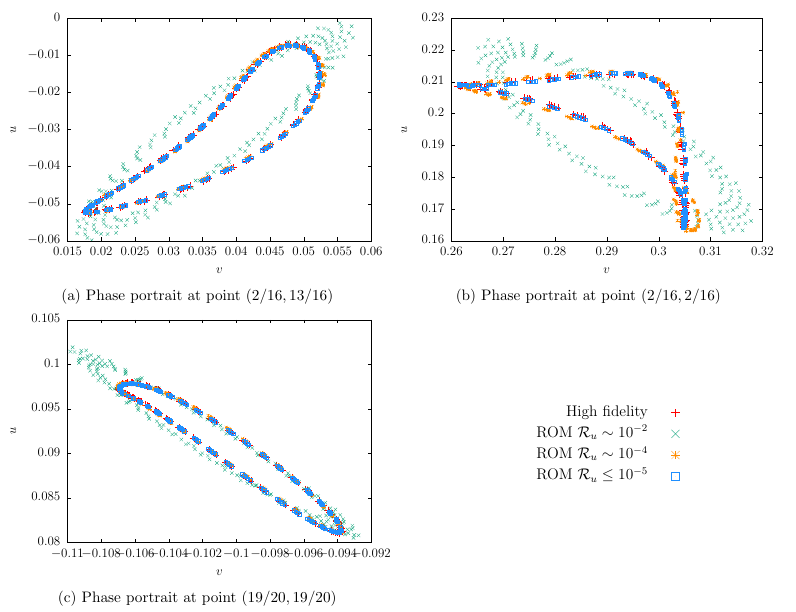}
\caption{Velocity phase portrait at the corner of the cavity domain for $Re = 9500$. $u$ and $v$ are respectively the horizontal and vertical components of the velocity vector $\bm{u}$.}
\label{fig:velocity_at_points}
\end{figure}

The dynamics of the flow solutions is shown in figure \ref{Fig:flow cavity}. We observe that in the lower left corner, the secondary vortex separates into two small vortexes that reincorporate in a periodic manner. In the same way the secondary vortex in the upper left corner has a similar behavior. In order to construct the POD bases, \ModifMourad{$200$ snapshots uniformly selected from the periodic regime with $10$ periods are used}. In Figure \ref{RC_RE9500}, we show the evolution of the contribution ratio of POD modes with respect to modes numbers. 
The effect of the POD contribution ratio on the SM-ROM solution is illustrated in Figure \ref{fig:velocity_at_points}. It shows the phase portraits\footnote{The evolution of the horizontal velocity in terms of the vertical velocity} at points located at the bottom left $(2/16, 2/16)$, top left $(2/16, 13/16)$ and top right $(0.95, 0.95)$ of the cavity. Due to the periodic nature of the flow at these points, the curves of the phase portraits are closed and follow the same path over the periods of the flow.

Starting from a POD contribution ratio of order $10^{-2}$, it can be seen that the SM-ROM phase portraits are not periodic and do not match the high fidelity solution. By decreasing the ratio to $10^{-4}$, the SM-ROM phase portraits tend to converge to the high fidelity solution except for the the bottom left point of the cavity where small differences can still be noticed. These differences disappear with ratios  of order $10^{-5}$ or less.
This shows that the predictions by the SM-ROM converge to the groundtruth solution with decreasing values of the POD contribution ratio. 

\begin{figure}[hbtp!]
\centering 
\includegraphics[width=\linewidth]{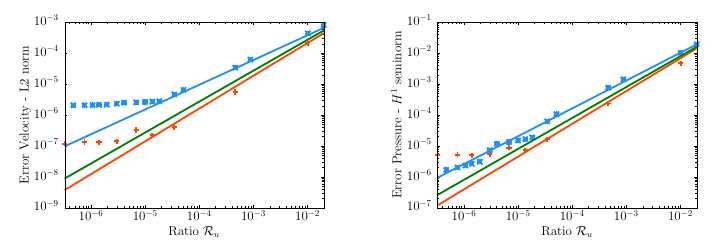}
\caption{Log-Log scale representation of the estimator versus SM-ROM errors for the flow in a lid driven cavity at $Re=9500$. The dots refer to the calculated values while the solid lines are their power regression. The green solid line corresponds to the curve of logarithmic slop equal to $1$. \ModifMourad{Only the errors corresponding to $\mcal{R}_u$ larger than $10^{-5}$ have been taken into account to calculate the regression lines.}}
\label{Errors_RE9500}
\end{figure}

In terms of error estimates, we observe from Figure \ref{Errors_RE9500}, that the estimator succeeds to predict the SM-ROM error slopes with a sufficiently good accuracy. 
%
The logarithmic slopes for the SM-ROM errors are close to one which gives an approximate rate of convergence of $10^{-1}$. 
%
\ModifMourad{
This is in good agreement with Figure \ref{fig:velocity_at_points} since between POD ratios $10^{-2}$ and $10^{-5}$, the SM-ROM jumps three orders of magnitudes and attains $10^{-5}$ for velocity and $4\times 10^{-6}$ for pressure. These errors represent a sufficient accuracy to reproduce well the dynamics of the flow and thus meet the high fidelity model hydrodynamics coefficients. 
}
As for the flow past a cylinder, the slopes of the estimators are smaller than the SM-ROM slopes. Nevertheless, they provide a good indicator of the convergence rate of the SM-ROM. Typically, at ratio $10^{-5}$, the estimator indicates an error value of $\frac{3}{2}\times 10^{-6}$ for velocity and $2\times 10^{-5}$ for the pressure. These are sufficiently good approximations of the SM-ROM velocity and pressure errors of values $\frac{3}{2}\times 10^{-7}$ and $4\times 10^{-6}$ respectively. Which confirms once again the robustness of the proposed estimator.


\section{Summary and conclusions}\label{sec:Concl}
\ModifMourad{
In this paper we have carried on the numerical analysis of the stabilized motivated reduced order (SM-ROM) solution of the pressure in incompressible flows. We treat the pressure equation by duality of the momentum conservation equation with gradients of the reduced pressure space, as function tests. We have used the Proper Orthogonal decomposition (POD) method with $L^2$ norms to build the reduced spaces for both velocity and pressure, and proven the stability of the velocity-pressure segregated SM-ROM discretization in adapted norms. In particular for pressure, the stability is proven in a discrete norm that has the same asymptotic order as the $L^2$ norm with respect to the mesh size. Both estimates for velocity and pressure are optimal with respect to the projection error on the POD reduced spaces. To validate the theoretical results, we have performed two two-dimensional numerical tests (lid-driven cavity flow and for flow past a cylinder), for which we nearly recover the optimality of the error estimates in  both problems. We observe an error stagnation effect for high ROM-POD dimensions, that possibly arises because the weak high-fidelity model is not exact for gradients of discrete pressures. Besides previous numerical investigation of the SM-ROM, the present analysis supports the use of this method to recover the pressure from ROMs of incompressible fluids that only retain weakly divergence-free velocities, with accuracy enough for practical applications of interest.}
\medskip

{\sl Acknowledgments:} This work has been partially supported by the Spanish Government-EU Feder grant RTI2018-093521-B-C31. T. Chac\'on deeply thanks the University of La Rochelle for its invitation in October of 2019 that allowed to start this research. The research of Samuele Rubino has been also funded by the Spanish State Research Agency through the national program Juan de la Cierva-Incorporaci\'on 2017.

\bibliographystyle{abbrv}
\bibliography{Biblio_SM-ROM}


\end{document}

%% file: NEW_COMMANDS.tex
\newcommand{\myfrac}[2]{\displaystyle{\frac{#1}{#2}}}
\newcommand{\somme}[3]{\displaystyle{\overset{#3}{\underset{#1=#2}{\sum}}}}

\newcommand{\vect}[1]{ {\bm{#1}}}

\newcommand{\mcal}[1]{\mathcal{#1}}
\newcommand{\mbb}[1]{\mathbb{#1}}

\newcommand{\txt}[1]{\textnormal{#1}}

%% file: SM-ROM.bbl
\begin{thebibliography}{10}

\bibitem{metstab}
N.~Ahmed, T.~Chac\'{o}n~Rebollo, V.~John, and S.~Rubino.
\newblock Analysis of a full space-time discretization of the {N}avier-{S}tokes
  equations by a local projection stabilization method.
\newblock {\em IMA J. Numer. Anal.}, 37(3):1437--1467, 2017.

\bibitem{ALDMOUR2002}
A.~AL-Dmour and K.~Mohammad.
\newblock Active control of flexible structures using principal component
  analysis in the time domain.
\newblock {\em Journal of Sound and Vibration}, 253, 2002.

\bibitem{Allery2008}
C.~Allery, C.~Beghein, and A.~Hamdouni.
\newblock On investigation of particle dispersion by a \textsc{POD} approach.
\newblock {\em International Applied Mechanics}, 44(1):110--119, 2008, \textit{
  https$:$//doi.org/10.1007/s10778-008-0025-2}.

\bibitem{Allery-2005}
C.~Allery, C.~Béghein, and A.~Hamdouni.
\newblock Applying {P}roper {o}rthogonal {d}ecomposition to the computation of
  particle dispersion in a two-dimensional ventilated cavity.
\newblock {\em Communications in Nonlinear Science and Numerical Simulation},
  10:907--920, 2005.

\bibitem{Allery2005}
C.~Allery, C.~Béghein, and A.~Hamdouni.
\newblock Applying \textsc{P}roper \textsc{O}rthogonal \textsc{D}ecomposition
  to the computation of particle dispersion in a two-dimensional ventilated
  cavity.
\newblock {\em Communications in Nonlinear Science and Numerical Simulation},
  10(8):907--920, 2005, \textit{
  https$:$//doi.org/10.1016/j.cnsns.2004.05.005}.

\bibitem{Ballarin2020}
F.~Ballarin, T.~Chac\'{o}n~Rebollo, E.~Delgado~\'{A}vila,
  M.~G\'{o}mez~M\'{a}rmol, and G.~Rozza.
\newblock Certified reduced basis {VMS}-{S}magorinsky model for natural
  convection flow in a cavity with variable height.
\newblock {\em Comput. Math. Appl.}, 80(5):973--989, 2020.

\bibitem{Quarter2015}
F.~Ballarin, A.~Manzoni, A.~Quarteroni, and G.~Rozza.
\newblock Supremizer stabilization of podgalerkin approximation of parametrized
  steady incompressible navier-stokes equations.
\newblock {\em Internat. J. Numer. Methods Engrg.}, 102:1136--1161, 2015.

\bibitem{Berkooz}
G.~Berkooz, P.~Holmes, and J.~L. Lumley.
\newblock On the relation between low-dimensional models and the dynamics of
  coherent structures in the turbulent wall layer.
\newblock {\em Theor. Comput. Fluid Dyn.}, 4:255--269, 1993.

\bibitem{Bernardi04}
C.~Bernardi, Y.~Maday, and F.~Rapetti.
\newblock {\em Discr{\'e}tisations variationnelles de probl{\`e}mes aux limites
  elliptiques}, volume~45 of {\em Math{\'e}matiques \& Applications}.
\newblock Springer-Verlag, 2004.

\bibitem{Brezis11}
H.~Brezis.
\newblock {\em Functional analysis, {S}obolev spaces and partial differential
  equations}.
\newblock Universitext. Springer, New York, 2011.

\bibitem{BF}
F.~Brezzi and R.~S. Falk.
\newblock Stability of higher-order {H}ood-{T}aylor methods.
\newblock {\em SIAM J. Numer. Anal.}, 28(3):581--590, 1991.

\bibitem{Buljak2011}
V.~Buljak and G.~Maier.
\newblock {P}roper {o}rthogonal {d}ecomposition and radial basis functions in
  material characterization based on instrumented indentation.
\newblock {\em Engineering Structures}, 33, 2011.

\bibitem{Beghein2014}
C.~Béghein, C.~Allery, J.~Pozorski, and M.~Waclawczyk.
\newblock Application of \textsc{POD}-based dynamical systems to dispersion and
  deposition of particles in turbulent channel flow.
\newblock {\em International Journal of Multiphase Flow}, 58:97 -- 113, 2014,
  \textit{ https$:$//doi.org/10.1016/j.ijmultiphaseflow.2013.09.001}.

\bibitem{Beghein-Allery-2014}
C.~Béghein, C.~Allery, M.~Wacławczyk, and J.~Pozorski.
\newblock Application of {POD}-based dynamical systems to dispersion and
  deposition of particles in turbulent channel flow.
\newblock {\em International Journal of Multiphase Flow}, 58:97--113, 1 2014.

\bibitem{IliescuJohn14}
A.~Caiazzo, T.~Iliescu, V.~John, and S.~Schyschlowa.
\newblock A numerical investigation of velocity-pressure reduced order models
  for incompressible flows.
\newblock {\em J. Comput. Phys.}, 259:598--616, 2014.

\bibitem{Cha_maca}
T.~Chac{\'o}n~Rebollo, M.~G{\'o}mez~M{\'a}rmol, V.~Girault, and
  I.~S{\'a}nchez~Mu{\~n}oz.
\newblock A high order term-by-term stabilization solver for incompressible
  flow problems.
\newblock {\em IMA J. Numer. Anal.}, 33(3):974--1007, 2013.

\bibitem{chalew}
T.~Chac\'{o}n~Rebollo and R.~Lewandowski.
\newblock {\em Mathematical and numerical foundations of turbulence models and
  applications}.
\newblock Modeling and Simulation in Science, Engineering and Technology.
  Birkh\"{a}user/Springer, New York, 2014.

\bibitem{Ciarlet78}
P.~Ciarlet.
\newblock {\em The finite element method for elliptic problems.}
\newblock Studies in Mathematics and its Applications, Vol. 4. North-Holland
  Publishing Co., Amsterdam-New York-Oxford, 1978.

\bibitem{NS_grad_div}
J.~de~Frutos, B.~Garc\'{\i}a-Archilla, V.~John, and J.~Novo.
\newblock Analysis of the grad-div stabilization for the time-dependent
  {N}avier-{S}tokes equations with inf-sup stable finite elements.
\newblock {\em Adv. Comput. Math.}, 44(1):195--225, 2018.

\bibitem{BoscoJulia19}
J.~de~Frutos, B.~Garc\'{\i}a-Archilla, and J.~Novo.
\newblock Fully discrete approximations to the time-dependent {N}avier-{S}tokes
  equations with a projection method in time and grad-div stabilization.
\newblock {\em J. Sci. Comput.}, 80(2):1330--1368, 2019.

\bibitem{BoscoJulia19Corr}
J.~de~Frutos, B.~Garc\'{\i}a-Archilla, and J.~Novo.
\newblock Corrigenda: {F}ully discrete approximations to the time-dependent
  {N}avier-{S}tokes equations with a projection method in time and grad-div
  stabilization.
\newblock {\em J. Sci. Comput.}, 88(2):Paper No. 40, 3, 2021.

\bibitem{Fick2018}
L.~Fick, Y.~Maday, A.~T. Patera, and T.~Taddei.
\newblock A stabilized {POD} model for turbulent flows over a range of
  {R}eynolds numbers: optimal parameter sampling and constrained projection.
\newblock {\em J. Comput. Phys.}, 371:214--243, 2018.

\bibitem{IliescuJohn15}
S.~Giere, T.~Iliescu, V.~John, and D.~Wells.
\newblock S{UPG} reduced order models for convection-dominated
  convection-diffusion-reaction equations.
\newblock {\em Comput. Methods Appl. Mech. Engrg.}, 289:454--474, 2015.

\bibitem{Han2003}
S.~Han and B.~Feeny.
\newblock Application of {P}roper {o}rthogonal {d}ecomposition to structural
  vibration analysis.
\newblock {\em Mechanical Systems and Signal Processing}, 17, 2003.

\bibitem{Hesthaven2016}
J.~S. Hesthaven, G.~Rozza, and B.~Stamm.
\newblock {\em Certified reduced basis methods for parametrized partial
  differential equations}.
\newblock SpringerBriefs in Mathematics. Springer, Cham; BCAM Basque Center for
  Applied Mathematics, Bilbao, 2016.
\newblock BCAM SpringerBriefs.

\bibitem{Hijazi2020}
S.~Hijazi, G.~Stabile, A.~Mola, and G.~Rozza.
\newblock Data-driven {POD}-{G}alerkin reduced order model for turbulent flows.
\newblock {\em J. Comput. Phys.}, 416:109513, 30, 2020.

\bibitem{Holmes96}
P.~Holmes, J.~L. Lumley, and G.~Berkooz.
\newblock {\em Turbulence, coherent structures, dynamical systems and
  symmetry}.
\newblock Cambridge Monographs on Mechanics. Cambridge University Press,
  Cambridge, 1996.

\bibitem{IliescuWang14}
T.~Iliescu and Z.~Wang.
\newblock Variational multiscale proper orthogonal decomposition:
  {N}avier-{S}tokes equations.
\newblock {\em Numer. Methods Partial Differential Equations}, 30(2):641--663,
  2014.

\bibitem{Kean2020}
K.~Kean and M.~Schneier.
\newblock Error analysis of supremizer pressure recovery for pod based
  reduced-order models of the time dependent navier-stokes equations.
\newblock {\em SIAM J. Numer. Anal.}, 58 (4):2235--2264, 2020.

\bibitem{KunischVolkwein01}
K.~Kunisch and S.~Volkwein.
\newblock Galerkin proper orthogonal decomposition methods for parabolic
  problems.
\newblock {\em Numer. Math.}, 90(1):117--148, 2001.

\bibitem{FenicsBook}
H.~P. Langtangen and A.~Logg.
\newblock {\em Solving PDEs in Python: The FEniCS Tutorial I}.
\newblock Springer Publishing Company, Incorporated, 1st edition, 2017.

\bibitem{Leblond2011}
C.~Leblond, C.~Allery, and C.~Inard.
\newblock An optimal projection method for the reduced-order modeling of
  incompressible flows.
\newblock {\em Computer Methods in App. Mechanics and Engineering},
  200(33-36):2507--2527, 2011, \textit{
  https$:$//doi.org/10.1016/j.cma.2011.04.020}.

\bibitem{Lucia2004}
D.~J. Lucia, P.~S. Beran, and W.~A. Silva.
\newblock Reduced-order modeling: new approaches for computational physics.
\newblock {\em Progress in Aerospace Sciences}, 40, 2004.

\bibitem{NovoRubinoSINUM21}
J.~Novo and S.~Rubino.
\newblock Error analysis of proper orthogonal decomposition stabilized methods
  for incompressible flows.
\newblock {\em SIAM J. Numer. Anal.}, 59(1):334--369, 2021.

\bibitem{Park2007}
S.~Park, J.-J. Lee, C.-B. Yun, and D.~J. Inman.
\newblock Electro-mechanical impedance-based wireless structural health
  monitoring using {PCA}-{D}ata compression and {k}-{m}eans {C}lustering
  {A}lgorithms.
\newblock {\em Journal of Intelligent Material Systems and Structures}, 19, 05
  2007.

\bibitem{AlfioRozza2009}
A.~Quarteroni and G.~Rozza.
\newblock Numerical solution of parametrized navier–stokes equations by
  reduced basis methods.
\newblock {\em Numerical Methods for Partial Differential Equations},
  23(4):923--948, 2007.

\bibitem{Rowley}
C.~W. Rowley, T.~Colonius, and R.~M. Murray.
\newblock Model reduction for compressible flows using pod and galerkin
  projection.
\newblock {\em Phys. D: Nonlinear Phenom.}, 189 (1):115--129, 2004.

\bibitem{Veroy2007}
G.~Rozza and K.~Veroy.
\newblock On the stability of the reduced basis method for stokes equations in
  parametrized domains.
\newblock {\em Comput. Methods Appl. Mech. Engrg.}, 196:1244--1260, 2007.

\bibitem{Rubino20}
S.~Rubino.
\newblock Numerical analysis of a projection-based stabilized {POD}-{ROM} for
  incompressible flows.
\newblock {\em SIAM J. Numer. Anal.}, 58(4):2019--2058, 2020.

\bibitem{Sirovich}
L.~Sirovich.
\newblock Turbulence and the dynamics of coherent structures : {P}art {I}, {II}
  and {III}.
\newblock {\em Quarterly of Applied Mathematics}, pages 461--590, 1987.

\bibitem{Tallet2016}
A.~Tallet, C.~Allery, and C.~Leblond.
\newblock Optimal flow control using a \textsc{POD}-based reduced-order model.
\newblock {\em Numerical Heat Transfer, Part B: Fundamentals}, 70(1):1--24,
  2016, \textit{ https$:$//doi.org/10.1080/10407790.2016.1173472}.

\bibitem{Tallet2015a}
A.~Tallet, C.~Allery, C.~Leblond, and E.~Liberge.
\newblock A minimum residual projection to build coupled velocity pressure
  \textsc{POD-ROM} for incompressible \textsc{N}avier-\textsc{S}tokes
  equations.
\newblock {\em Communications in Nonlinear Science and Numerical Simulation},
  22(1-3):909--932, 2015, \textit{
  https$:$//doi.org/10.1016/j.cnsns.2014.09.009}.

\bibitem{hood0}
C.~Taylor and P.~Hood.
\newblock A numerical solution of the {N}avier-{S}tokes equations using the
  finite element technique.
\newblock {\em Internat. J. Comput. \& Fluids}, 1(1):73--100, 1973.

\bibitem{Thomas2003}
J.~P. Thomas, E.~H. Dowell, and K.~C. Hall.
\newblock Three-dimensional transonic aeroelasticity using {P}roper
  {o}rthogonal {d}ecomposition-based reduced-order models.
\newblock {\em Journal of Aircraft}, 40, 05 2003.

\end{thebibliography}
